\documentclass{article}
\usepackage[utf8]{inputenc}
\usepackage{authblk}
\usepackage{setspace}
\usepackage[margin=1in]{geometry}
\usepackage{graphicx}
\graphicspath{ {./figures/} }
\usepackage{subcaption}
\usepackage{amsmath}
\usepackage{lineno}
\usepackage{comment}
\usepackage{amsthm}

\newtheorem{thm}{Theorem}[section]

\newtheorem{lem}[thm]{Lemma}

\usepackage{amsfonts} 
\usepackage{mathtools} 

\usepackage{url}
\usepackage{hyperref}
\usepackage{xcolor}
\usepackage{float}

\hypersetup{
	colorlinks,
	linkcolor={red!50!black},
	citecolor={blue!50!black},
	urlcolor={blue!80!black}
}
\usepackage[natbib=true,backend=biber,sorting=none,style=ieee, citestyle=numeric-comp]{biblatex}

\title{Global Dynamics and Time-Optimal Control Studies for Additional Food provided Holling Type-III Mutually Interfering Prey-Predator Systems with Applications to Pest Management}

\author[1]{D. Bhanu Prakash}
\author[2]{D. K. K. Vamsi}

\affil[1, 2]{ \ Department of Mathematics and Computer Science, Sri Sathya Sai Institute of Higher Learning, India.}
\affil[1]{First Author. Email: dbhanuprakash@sssihl.edu.in}
\affil[2]{Corresponding author. Email: dkkvamsi@sssihl.edu.in}

\date{}

\begin{document}
	
	\maketitle
	
	\begin{abstract} {{
				\noindent 
				In this study, we investigate a prey-predator model exhibiting Holling type-III functional response among mutually interfering predators to assess the effects of provision of additional food to natural enemies in altering pest dynamics. We derive and study both the non-additional food provided system (initial system) and additional food provided system mathematically and through extensive numerical simulations. We prove the positivity and boundedness of the solutions of both the systems. We perform an extensive existence and stability analysis of various equilibrium points exhibited by both the systems. We discuss the global dynamics and stability behavior exhibited by both the systems with respect to crucial parameters using numerical methods. Considering the quality and quantity of additional food as control parameters, we framed and characterized the time optimal control problems. We further simulated these time optimal problems using numerical optimization techniques. Our theoretical and numerical investigations reveal that the provision of suitable choice of additional food steers the system to a prey-elimination state, leading to a pest-regulated ecosystems.
		}}
	\end{abstract}
	
	{ \bf {Article highlights:} }
	
	\begin{itemize}
		\item The impact of the provision of additional food to the Predators (Natural Enemies) of Prey (Pests) was assessed using theoretical and numerical analysis of a Holling Type III mathematical model. 
		\item Additional food and the mutual interference among predators was observed to significantly influence the populations of prey, predator over time.
		\item The optimal quality(quantity) of additional food required to reach pest-regulated state in least possible time is observed to be only in the extremes in the possible interval. 
	\end{itemize}
	

	
	
	\section{Introduction}
	
	For several decades, the control of economically harmful agricultural weeds and insect pests has garnered significant attention from researchers across multiple disciplines. It is widely recognized that the use of pesticides and insecticides can have enduring negative impacts on both the environment and human health. Alternative pest management techniques such as the deployment of natural enemies of pests are widely researched in the past few decades. Provision of additional food for mass rearing and releasing the natural enemies into the fields or providing them as supplementary food in the fields have proven to be efficient in controlling the pest populations \cite{harwood2004prey,putman2004supplementary,redpath2001does,winkler2005plutella}. 
	
	Mathematical modeling of the prey-predator systems have significantly contributed to the enhancement of understanding of this strategy. Natural enemies are different for different types of pests unlike many pesticides which can be used for multiple pests in general. The functional response i.e., the behavior at which natural enemy feeds on pest, plays a crucial role in mathematical modeling of the pest-natural enemy interactions. The impact of additional food on natural enemies that follow Holling type-II functional response are well studied in literature \cite{varaprasadsiraddi,varaprasadsirquant,varaprasadsirtime}. Recently, authors in \cite{V3EarthSystems,V3JTB,V4Acta,V4DEDS} studied the systems exhibiting Holling type-III and Holling type-IV functional responses. 
	
	However, most of these models reach infinite prey in the absence of predator. Also, these models does not incorporate the intra-specific competition among predators. As the ecosystems are of limited resources, competition plays an important role in the dynamics of the system. This mutual interference can be incorporated in two ways. Bazykin models incorporate this feature explicitly \cite{bazykin1976structural}. In this paper, we followed the method from \cite{varaprasadsiraddi} and incorporated the mutual interference implicitly for the model exhibiting Holling type-III functional response. 
	
	This paper is structured as follows: In the next section, the Holling type-III functional response among mutual interfering predators is derived. In the following Section \ref{secderm}, the initial system without additional food and the prey-predator model in the presence of additional food to predators is developed. The positivity and boundedness of these systems is proved in section \ref{secposi}. In section \ref{secexis}, the predator-prey isoclines and the existence of equilibria are discussed for both the systems. In section \ref{secstab}, local stability of their equilibria are analyzed for these systems. Sections \ref{secglob1} and \ref{secglob2} discusses the global dynamics of initial system and additional food provided system respectively. Section \ref{secmutual} discusses the effect of mutual interference on the dynamics of the additional food provided system. Section \ref{seccontrol} formulates the time optimal control problem with quality and quantity of additional food as control parameters and characterizes the optimal controls. The results obtained here are numerically simulated in section \ref{secnumer}. Finally, in section \ref{secdisc} the discussion and conclusions are presented.

	\section{Derivation of Holling Type-III Functional Response} \label{secfr}
	Let $N,\ P$ denote the densities of prey and predator respectively. 
	
	Let $\triangle T$, a small period of time (small in the sense that the predator and prey densities remain roughly constant over $\triangle$T) that the predator spends for searching prey and/or additional food, consuming the captured prey and/or additional food, and interacting with other predators.
	
	Let $\triangle T_S$ denote the part of $\triangle T$ that the predator spends for searching prey and/or additional food. 
	
	Let $\triangle T_N$ denote the part of $\triangle T$ that the predator spends for handling the prey.
	
	Let $\triangle T_A$ denote the part of $\triangle T$ that the predator spends for handling the additional food. 
	
	Let $\triangle T_P$ denote the part of $\triangle T$ that the predator spends for interacting with other predators.
	
	So, we have
	\begin{equation} \label{delt3}
		\triangle T = \triangle T_S + \triangle T_N + \triangle T_A + \triangle T_P.
	\end{equation}
	
	Let $h_N, h_A$ and $h_P$ denote the length of time required for each interaction between a predator and prey item, additional food, other predators respectively.
	
	Let $e_N$ represents the search rate of the predator per unit prey availability and $e_A$ represents the search rate of the predator per unit quantity of additional food. Let $e_P$ represents the rate constant at which a predator encounters other predators.
	
	From \cite{V3JTB}, the total handling time for the prey ($\triangle T_N$) equals $h_N$ times number of prey caught, which is given by, $h_N e_N N^2 \triangle T_S$. 
	
	Now, the handling time for the additional food equals the handling time for one additional food item times the total density of additional food encountered. Also,the additional food encountered is proportional to the search time and the additional food density. Hence, the quantity of additional food encountered during this time is proportional to $\triangle T_S  A = e_A A \triangle T_S$, where $e_A$ is the proportionality constant, which denotes the catchability of the additional food. Now, the total handling time for the additional food ($\triangle T_A$)  equals $h_A$ times additional food encountered, which is given by, $h_A e_A A \triangle T_S$. 
	
	Similarly, $\triangle T_P$ is given by $h_P e_P P  \triangle T_S$.
	
	Hence, 
	\begin{eqnarray}
		\triangle T_N &=& h_N e_N N^2 \triangle T_S \label{tn3} \\
		\triangle T_A &=& h_A e_A A  \triangle T_S \label{ta3} \\
		\triangle T_P &=& h_P e_P P  \triangle T_S \label{tp3} 
	\end{eqnarray}
	From (\ref{delt3}), (\ref{tn3}), (\ref{ta3}) and (\ref{tp3}), we have
	
	\begin{equation} \label{deltt3}
		\triangle T = \Big(1 + h_N e_N N^2 + h_A e_A A  + h_P e_P P \Big) \triangle T_S
	\end{equation}
	
	Now, as a predator encounters $e_N N^2 \triangle T_S$ prey items during the period $\triangle T$, the overall rate of encounters with prey over the time interval $\triangle T$ is given by
	
	\begin{equation*}
		\begin{split}
			g(N,P,A) & = \frac{\text{Total number of Prey Caught}}{\text{Total Time}}\\
			& = \frac{e_N N^2 \triangle T_S}{\triangle T} = \frac{e_N N^2 \triangle T_S}{(1 + h_N e_N N^2  + h_A e_A A  + h_P e_P P ) \triangle T_S} \\
			& = \frac{e_N N^2}{1 + h_N e_N N^2 + h_A e_A A  + h_P e_P P} \\
			& \text{Dividing numerator and denominator by } e_N h_N \\
			& = \frac{\frac{1}{h_N} N^2}{\frac{1}{e_N h_N} + N^2 + \frac{h_A e_A}{h_N e_N} A  + \frac{h_P e_P}{h_N e_N} P} \\
		\end{split}
	\end{equation*}
	
	Similarly, as a predator encounters $e_A A \triangle T_S$ quantity of additional food during the period $\triangle T$, the overall rate of encounters with additional food over the time interval $\triangle T$ is given by
	
	\begin{equation*}
		\begin{split}
			h(N,P,A) & = \frac{\text{Total Additional Food Encountered}}{\text{Total Time}}\\
			& = \frac{e_A A \triangle T_S}{\triangle T} = \frac{e_A A \triangle T_S}{(1 + h_N e_N N^2 + h_A e_A A  + h_P e_P P) \triangle T_S} \\
			& = \frac{ e_A A}{1 + h_N e_N N^2 + h_A e_A A  + h_P e_P P} \\
			& \text{Dividing numerator and denominator by } e_N h_N \\
			& = \frac{\frac{e_A }{e_N h_N} A }{\frac{1}{e_N h_N}+ N^2 + \frac{h_A e_A}{h_N e_N}A+\frac{h_P e_P}{h_N e_N} P}\\
		\end{split}
	\end{equation*}
	
	Hence, the Holling type-III functional response of the mutually interfereing predator towards the prey ($g(N,P,A)$) and the additional food ($h(N,P,A)$) respectively are given by
	
	\begin{eqnarray}
		g(N,P,A) = \frac{\frac{1}{h_N} N^2}{\frac{1}{e_N h_N} + N^2 + \frac{h_A e_A}{h_N e_N} A  + \frac{h_P e_P}{h_N e_N} P} \label{g3} \\
		h(N,P,A) =  \frac{\frac{e_A }{e_N h_N} A }{\frac{1}{e_N h_N}+ N^2 + \frac{h_A e_A}{h_N e_N}A+\frac{h_P e_P}{h_N e_N} P} \label{h3}
	\end{eqnarray}
	
	\section{Derivation of Model} \label{secderm}
	
	Now, we derive the additional food provided mutually interfering prey-predator model representing the prey-predator dynamics when the predator is provided with additional food (which is of non-reproducing item) based on the functional response derived in the previous section. 
	
	The deterministic prey-predator model exhibiting Holling type-III functional response among mutually interfering predators and where prey grow logistically in the absence of predator is given by
	
	\begin{equation*}
		\begin{split}
			\frac{\mathrm{d}N}{\mathrm{d}T} & = r N \Big(1-\frac{N}{K}\Big) - g(P,N,A=0) \ P \\
			\frac{\mathrm{d}P}{\mathrm{d}T} & = \epsilon_N \ g(P,N,A=0)\  P - m_1 P\\
		\end{split}
	\end{equation*}
	
	Here $\epsilon_N$ is the conversion factor that represents the rate at which the prey biomass gets converted into predator biomass.  The parameters $r$ and $K$ represent the intrinsic growth rate and carrying capacity of the prey respectively. $m_1$ is the death rate of predator in the absence of prey which is also termed as starvation rate. The biological descriptions of the various parameters involved in the system (\ref{midms03}) is described in table \ref{param_tab3}.
	
	Substituting (\ref{g3}) and (\ref{h3}) in the above equations, we get 
	
	\begin{equation*}
		\begin{split}
			\frac{\mathrm{d}N}{\mathrm{d}T} & = r N \Big(1-\frac{N}{K}\Big) -  \frac{\frac{1}{h_N} N^2}{\frac{1}{e_N h_N} + N^2 +  \frac{h_P e_P}{h_N e_N} P} P \\
			\frac{\mathrm{d}P}{\mathrm{d}T} & = \Big( \epsilon_N  \frac{\frac{1}{h_N} N^2}{\frac{1}{e_N h_N} + N^2 + \frac{h_P e_P}{h_N e_N} P} \Big) P - m_1 P\\
		\end{split}
	\end{equation*}
	
	Here, let $c = \frac{1}{h_N}$ and $a = \frac{1}{\sqrt{h_N e_N}}$ stands for the maximum rate of predation and half-saturation value of the predator respectively. 
	
	So, by substituting $\delta_1 = \epsilon_N c$ and $\epsilon_1 = e_P h_P$, the model gets transformed to 
	
	\begin{equation}\label{midms03}
		\begin{split}
			\frac{\mathrm{d}N}{\mathrm{d}T} & = r N \Big(1-\frac{N}{K}\Big) -  \frac{cN^2P}{a^2+N^2+ \epsilon_1 a^2 P} \\
			\frac{\mathrm{d}P}{\mathrm{d}T} & = \delta_1 \Big( \frac{N^2}{a^2+N^2+ \epsilon_1 a^2 P} \Big)P - m_1 P \\
		\end{split}
	\end{equation}
	
	In similar lines to the above derivation, the deterministic prey-predator model exhibiting Holling type-III functional response in the presence of additional food for mutually interfering predators and where prey grow logistically in the absence of predator is given by
	
	\begin{equation*}
		\begin{split}
			\frac{\mathrm{d}N}{\mathrm{d}T} & = r N \Big(1-\frac{N}{K}\Big) - g(P,N,A) P \\
			\frac{\mathrm{d}P}{\mathrm{d}T} & = (\epsilon_N g(P,N,A)+\epsilon_A h(N,P,A))P - m_1 P\\
		\end{split}
	\end{equation*}
	
	Here $\epsilon_N$ and $\epsilon_A$ are conversion factors that represents the rate at which the prey/additional food biomass gets converted into predator biomass.
	
	Substituting (\ref{g3}) and (\ref{h3}) in the above equations, we get 
	
	\begin{equation*}
		\begin{split}
			\frac{\mathrm{d}N}{\mathrm{d}T} & = r N \Big(1-\frac{N}{K}\Big) -  \frac{\frac{1}{h_N} N^2}{\frac{1}{e_N h_N} + N^2 + \frac{h_A e_A}{h_N e_N} A  + \frac{h_P e_P}{h_N e_N} P} P \\
			\frac{\mathrm{d}P}{\mathrm{d}T} & = \Big( \epsilon_N  \frac{\frac{1}{h_N} N^2}{\frac{1}{e_N h_N} + N^2 + \frac{h_A e_A}{h_N e_N} A  + \frac{h_P e_P}{h_N e_N} P} \\ &  +\epsilon_A   \frac{\frac{e_A }{e_N h_N} A }{\frac{1}{e_N h_N}+ N^2 + \frac{h_A e_A}{h_N e_N}A+\frac{h_P e_P}{h_N e_N} P} \Big)P - m_1 P\\
		\end{split}
	\end{equation*}
	
	Here, let $c = \frac{1}{h_N}$ and $a = \frac{1}{\sqrt{h_N e_N}}$ stands for the maximum rate of predation and half-saturation value of the predator respectively. 
	
	The term $\alpha = \frac{\epsilon_N/h_N}{\epsilon_A/h_A}$, which is the ratio between the maximum growth rates of the predator when it consumes the prey and additional food respectively, indicates the relative efficiency of the predator to convert either of the available food into predator biomass. The value $\alpha$ can be seen to be an equivalent of \textbf{quality} of additional food. Let $\eta = \frac{e_A \epsilon_A}{e_N \epsilon_N}$ and the term $\eta A$ effectual food level. Here the term $\epsilon_1 = e_P h_P$ represents the strength of mutual interference between predators. So, by substituting $\epsilon_A \frac{e_A}{e_N} = \epsilon_N \eta, \ \frac{h_A e_A}{h_N e_N} = \alpha \eta, \ \delta_1 = \epsilon_N c$ and $\epsilon_1 = e_P h_P$, the model gets transformed to 
	
	\begin{equation}\label{midms3}
		\begin{split}
			\frac{\mathrm{d}N}{\mathrm{d}T} & = r N \Big(1-\frac{N}{K}\Big) -  \frac{cN^2P}{a^2+N^2+\alpha \eta A+ \epsilon_1 a^2 P} \\
			\frac{\mathrm{d}P}{\mathrm{d}T} & = \delta_1 \Big( \frac{N^2+\eta A}{a^2+N^2+\alpha \eta A+ \epsilon_1 a^2 P} \Big)P - m_1 P \\
		\end{split}
	\end{equation}
	
	The biological descriptions of the various parameters involved in the systems (\ref{midms03}), (\ref{midms3}) is described in table \ref{param_tab3}.
	
	\begin{table}[bht!]
		\centering
		\begin{tabular}{ccc}
			\hline
			Parameter & Definition & Dimension \\  
			\hline
			T & Time & time\\ 
			N & Prey density & biomass \\
			P & Predator density & biomass \\
			A & Additional food & biomass \\
			r & Prey intrinsic growth rate & time$^{-1}$ \\
			K & Prey carrying capacity & biomass \\
			c & Maximum rate of predation & time$^{-1}$ \\
			$\delta_1$ & Maximum growth rate of predator & time$^{-1}$ \\
			$m_1$ & Predator mortality rate & time$^{-1}$ \\
			\hline
		\end{tabular}
		\caption{Description of variables and parameters present in the systems (\ref{midms03}), (\ref{midms3})}
		\label{param_tab3}
	\end{table}
	
	In order to reduce the complexity of the model, we non-dimensionalize the systems (\ref{midms03}), (\ref{midms3}) using the following transformations, 
	$$t=rT,\ N=ax, \  P=\frac{ary}{c}$$
	
	Accordingly, systems (\ref{midms03}), (\ref{midms3}) gets reduced to the following non-dimensionalised systems respectively :
	
	\begin{equation} \label{midm03}
		\begin{split}
			\frac{\mathrm{d} x}{\mathrm{d} t} & = x \Bigg(1-\frac{x}{\gamma} \Bigg)- \Bigg( \frac{x^2y}{1+x^2+\epsilon y} \Bigg) \\
			\frac{\mathrm{d} y}{\mathrm{d} t} & = \delta \Bigg( \frac{x^2}{1 + x^2 + \epsilon y} \Bigg) y - m y \\
		\end{split}
	\end{equation}
	
	\begin{equation} \label{midm3}
		\begin{split}
			\frac{\mathrm{d} x}{\mathrm{d} t} & = x \Bigg(1-\frac{x}{\gamma} \Bigg)- \Bigg( \frac{x^2y}{1+x^2+\alpha \xi+\epsilon y} \Bigg) \\
			\frac{\mathrm{d} y}{\mathrm{d} t} & = \delta \Bigg( \frac{x^2 + \xi}{1 + x^2 + \alpha \xi + \epsilon y} \Bigg) y - m y \\
		\end{split}
	\end{equation}
	
	where $$\gamma = \frac{K}{a}, \ \xi = \frac{\eta A}{a^2}, \  \epsilon = \frac{\epsilon_1 a r}{c}, \ \delta = \frac{\delta_1}{c},\ m_1 = c m.$$
	
	Here the term $\frac{\eta A}{N}$ denotes the quantity of additional food perceptible to the predator with respect to the prey relative to the nutritional value of prey to the additional food. Hence the term $\xi = \frac{\eta A}{a^2}$ can be seen to be an equivalent of {\textit {\bf{quantity}} }of additional food.
	
	For the sake of brevity, we call the system (\ref{midm03}) as initial system as we do not incorporate the additional food component. We call the system (\ref{midm3}) as additional food provided system. 
	
	\section{Positivity and Bondedness of Solutions} \label{secposi}
	
	\subsection{Positivity of the solutions}
	
	From the model system (\ref{midm3}), it is observed that
	$$ \frac{\mathrm{d} x}{\mathrm{d} t} \Big|_{x=0} \geq 0,\ \frac{\mathrm{d} y}{\mathrm{d} t}\Big|_{y=0} \geq 0$$
	
	The non-negativity of both rates is observed on the bounding planes ($x = 0, \ y=0$) within the non-negative region of real space. Consequently, if a solution initiates within the confines of this region, it will persist within it over time $t$. This phenomenon occurs due to the vector field consistently pointing inward on the bounding planes, as implied by the inequalities mentioned above.
	
	\subsection{Boundedness of the solutions}
	
	\begin{thm}
		Every solution of the system (\ref{midm3}) that starts within the positive quadrant of the state space remains bounded.
	\end{thm}\label{bound3}
	
	\begin{proof}
		We define $W = x + \frac{1}{\delta}y$. Now, for any $K > 0$, we consider,
		\begin{equation*}
			\begin{split}
				\frac{\mathrm{d} W}{\mathrm{d} t} + K W = & \frac{\mathrm{d} x}{\mathrm{d} t} + \frac{1}{\delta} \frac{\mathrm{d} y}{\mathrm{d} t} + K x + \frac{K}{\delta} y \\
				=  & x \Bigg(1-\frac{x}{\gamma} \Bigg)- \Bigg( \frac{x^2y}{1+x^2+\alpha \xi+\epsilon y} \Bigg) \\ & +\frac{1}{\delta}  \Bigg( \delta \Bigg( \frac{x^2 + \xi}{1 + x^2 + \alpha \xi + \epsilon y} \Bigg) y - m y \Bigg) + K x + \frac{K}{\delta} y \\
				=  & x -\frac{x^2}{\gamma} - \Bigg( \frac{x^2y}{1+x^2+\alpha \xi+\epsilon y} \Bigg) \\ & +  \Bigg( \frac{x^2 + \xi}{1 + x^2 + \alpha \xi + \epsilon y} \Bigg) y - \frac{m}{\delta} y + K x + \frac{K}{\delta} y \\
				=  & (1+K) x -\frac{x^2}{\gamma} +  \Bigg( \frac{\xi y}{1 + x^2 + \alpha \xi + \epsilon y} \Bigg) + \frac{K-m}{\delta} y \\
				\leq  & \frac{\gamma (1+K)^2}{4} +  \frac{\xi}{\epsilon } + \frac{K-m}{\delta} y \\
			\end{split}
		\end{equation*}
		For sufficiently small $K(<m)$, we get
		$$\frac{\mathrm{d} W}{\mathrm{d} t} + K W \leq M \Bigg(= \frac{\gamma (1+K)^2}{4} +  \frac{\xi}{\epsilon} \Bigg).$$
		From Gronwall's inequality, we get 
		$$0 \leq W(t) \leq \frac{M}{K} (1 - e^{-Kt}) + W(0) e^{-Kt}.$$
		Therefore, $0 < W(t) \leq \frac{M}{K}$ as $t \rightarrow \infty$.
	\end{proof}
	
	\section{Existence of Equilibria} \label{secexis}
	
	In this section, we discuss the existence of equilibrium points for the systems (\ref{midm03}) and (\ref{midm3}). We also discuss the behavior of isoclines based on the characteristics of the additional food. We consider the biologically feasible parametric constraint $\delta > m$. 
	
	The prey nullclines of the system (\ref{midm3}) are given by 
	$$ x = 0 , \  \ 1- \frac{x}{\gamma}  - \frac{x y}{1+x^2 + \alpha \xi + \epsilon y} = 0 $$
	
	The predator nullclines are given by
	$$ y = 0, \ \ \frac{\delta (x^2 + \xi)}{1+x^2+ \alpha \xi + \epsilon y} - m = 0$$
	
	Upon simplification, the non trivial predator nullcline is a parabola facing upwards $$y = \frac{(\delta - m) x^2 + \delta \xi - m (1+ \alpha \xi)}{m \epsilon}$$
	
	This parabola touches $y$-axis only at $\big(0,\frac{\delta \xi - m (1+\alpha \xi)}{m \epsilon}\big)$. In the absence of additional food, the nullcline touches negative $y$-axis i.e., $(0,\frac{-1}{\epsilon})$. With the provision of additional food, the predator nullcline moves upwards (when $\alpha < \frac{\delta}{m}$ and $0 < \xi < \frac{m}{\delta - m \alpha}$) or downwards (when $\alpha > \frac{\delta}{m}$ or $0<\frac{m}{\delta - m \alpha} < \xi$, if $\alpha < \frac{\delta}{m}$). It also passes through the positive $x$-axis at $(\sqrt{\frac{\delta \xi - m(1+\alpha \xi)}{m - \delta}},0)$ and touches the $x$-axis only when $\delta \xi - m(1+\alpha \xi) < 0$. 
	
	Upon simplification, the non trivial prey isocline is given by $$y = \frac{\Big(1-\frac{x}{\gamma}\Big)(1+x^2+\alpha \xi)}{\Big(1 + \frac{\epsilon}{\gamma}\Big) x - \epsilon }$$ 
	
	This prey nullcline passes through $(\gamma, 0)$ in positive quadrant. The line $x = \frac{\epsilon}{1+\frac{\epsilon}{\gamma}} $ becomes an asymptote for the prey isocline of (\ref{midm03}) and (\ref{midm3}). The prey nullcline lies in the positive $xy$-quadrant only when $\frac{\epsilon}{1+\frac{\epsilon}{\gamma}} < x < \gamma$ and remains negative for the remaining $x$ values. The prey nullcline has a slope
	
	\begin{eqnarray*} 
		\frac{\mathrm{d} y}{\mathrm{d} x} &=& \frac{\mathrm{d} }{\mathrm{d} x} \bigg(1 - \frac{x}{\gamma} \bigg) \frac{1+x^2 + \alpha \xi}{(1+\frac{\epsilon}{\gamma})x - \epsilon} +    \bigg(1 - \frac{x}{\gamma} \bigg) \frac{\mathrm{d} }{\mathrm{d} x} \bigg( \frac{1+x^2 + \alpha \xi}{(1+\frac{\epsilon}{\gamma})x - \epsilon} \bigg) \\
		&=& \frac{1}{\gamma \big((1+\frac{\epsilon}{\gamma})x - \epsilon\big)^2} \bigg[-2\bigg(1 + \frac{\epsilon}{\gamma}\bigg) x^3 + (4 \epsilon + \gamma) x^2 - 2 \epsilon \gamma x - \gamma (1+\alpha \xi) \bigg]
	\end{eqnarray*}
	
	The slope of prey isocline is a cubic equation in $x$. The discriminant of this cubic equation is given by $$\triangle = 4 (\gamma ^2 - 8 \epsilon \gamma - 27 (1+\alpha \xi)) (\epsilon^2 \gamma^2 + (\epsilon + \gamma)^2 (1+\alpha \xi)).$$
	
\begin{figure}[ht]
		\centering
		\includegraphics[width=\textwidth]{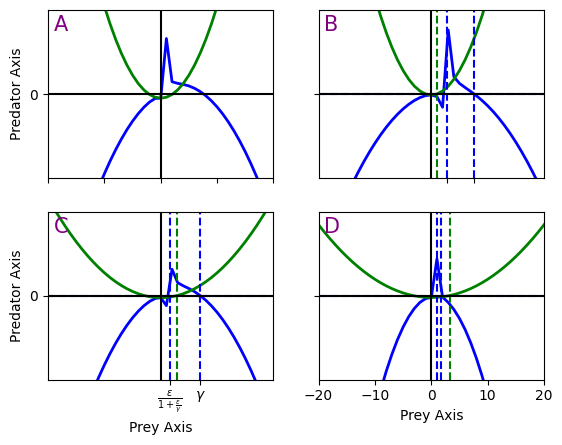}
		\caption{The possible configurations for the prey nullcline (represented by solid blue line) and predator nullcline (represented by solid red line) of the system (\ref{midm3}). The solid green line represents the predator nullcline of the initial system (\ref{midm03}).  }
		\label{iso13}
	\end{figure}

	The prey nullcline also exhibits a crest and trough as long as $\gamma > 4 \epsilon + \sqrt{16 \epsilon^2 + 27 (1+\alpha \xi)}$ (i.e., $\triangle > 0$) and it will be monotonically decreasing otherwise (i.e., $\triangle < 0$). 
	The possible configurations for the prey isocline and predator isocline are presented in figure \ref{iso13}.
	
	In figure \ref{iso13}, the vertical lines are $x = \frac{\epsilon}{1+\frac{\epsilon}{\gamma}}$ and $x = \gamma$. The frames A and C has predator nullcline of additional food provided system (\ref{midm3}) is below the predator nullcline of the initial system (\ref{midm03}). Whereas, the predator nullcline of additional food provided system (\ref{midm3}) is above that of the initial system (\ref{midm03}) in frames B and D. This implies that the interior equilibrium of the additional food provided system can appear only for $x > \frac{\delta \xi - m (1+ \alpha \xi)}{m - \delta} > 0$ in cases A and C, whereas in cases B and D, it can be any $x>0$.
	
	In the case of prey nullcline, frames A and B exhibit a crest and trough by satisfying the condition $\gamma > 4 \epsilon + \sqrt{16 \epsilon^2 + 27 (1+\alpha \xi)}$. Whereas, the frames C and D exhibits a steep decline to zero as x reaches $\gamma$. From the qualitative behaviour of the nullclines, it is observed that at least one interior equilibrium point exists in the positive-$xy$ quadrant for cases B and D. Also, interior equilibrium exists for cases A and C only when $\frac{\delta \xi - m (1+ \alpha \xi)}{m - \delta} < \gamma$.
	
	From the qualitative theory of dynamical systems, the interior equilibrium point exists only when $\sqrt{\frac{\delta \xi - m(1+\alpha \xi)}{m - \delta}} < \gamma$. Upon simplification, we get
	
	\begin{equation} \label{condiso3}
		\delta \xi - m(1+\alpha \xi) > - \gamma^2 (\delta - m)
	\end{equation}
	
	The systems (\ref{midm03}) and (\ref{midm3}) always admits the points $E_0 = (0,0)$ and $E_1 = (\gamma , 0)$ as their trivial and axial equilibrium respectively. The system  (\ref{midm3}) exhibits another axial equilibrium $E_2 = (0,\frac{\delta \xi - m (1+\alpha \xi)}{m \epsilon})$ if $\delta \xi - m(1+\alpha \xi) > 0$. The interior equilibrium of the system (\ref{midm03}), if exists, is given by $\bar{E} = (\bar{x}, \bar{y})$  where 
	
	$$ \bar{y} = \frac{(\delta - m) \bar{x}^2 - m}{m \epsilon}$$ 
	and 
	\begin{eqnarray*}
		& & (\epsilon \delta + \gamma (\delta - m)) \bar{x}^2 - \gamma \epsilon \delta \bar{x} - \gamma m = 0 \\
		\implies & & \bar{x} = \frac{\epsilon \gamma \delta + \sqrt{(\epsilon \gamma \delta)^2 + 4 m \gamma ((\delta - m) \gamma + \delta \epsilon)}}{2 ((\delta - m)\gamma + \delta \epsilon)} \  > 0 \\
	\end{eqnarray*}
	
	The interior equilibrium of the system (\ref{midm3}), if exists, is given by $E^* = (x^*, y^*)$  where 
	
	$$ y^* =  \frac{(\delta - m) (x^*)^2 + \delta \xi - m (1+\alpha \xi)}{m \epsilon}$$ 
	and $x^*$ should satisfy the following cubic equation
	$$ (\gamma (\delta  - m) + \epsilon \delta) (x^*)^3 - \epsilon \delta \gamma (x^*)^2 + (\epsilon \delta \xi + \gamma (\delta \xi - m (1+\alpha \xi))) x^* - \epsilon \delta \gamma \xi = 0.$$
	
	Comparing this with the standard cubic equation $a x^3 + b x^2 + c x + d = 0$, we have $a > 0,\  b<0,\ d<0$. This implies that the sum of roots of this cubic equation is positive and the product of roots of this cubic equation is also positive. From the theory of cubic equations, there exists at least one real root for any cubic equation. Since the product of roots is positive, there exists atleast one positive real root. This guarentees the existence of interior equilibrium point for (\ref{midm3}) provided $y^* > 0$.   
	
	We also have $\bar{x} < \gamma$ and $x^* < \gamma$. In the absence of additional food ($\alpha = \xi = 0$), the interior equilibrium points $\bar{E} = E^*$.
	
	To get more insights into the dynamics of these prey-predator models, we perform the stability analysis in the following section. 
	
	\section{Stability of Equilibria} \label{secstab}
	
	In order to obtain the asymptotic behavior of the trajectories of the system (\ref{midm3}), the associated Jacobian matrix is given by 
	
	$$J = \left[
	\begin{matrix}
		\frac{\partial}{\partial x} f(x,y)  & \frac{\partial}{\partial y} f(x,y)\\
		\frac{\partial}{\partial x} g(x,y) & \frac{\partial}{\partial y} g(x,y)
	\end{matrix}
	\right]$$
	
	where
	
	$$f(x,y) = x \Bigg(1-\frac{x}{\gamma} \Bigg)- \Bigg( \frac{x^2y}{1+x^2+\alpha \xi +\epsilon y} \Bigg) $$
	
	$$g(x,y) = \delta \Bigg( \frac{x^2+\xi}{1 + x^2 + \alpha \xi + \epsilon y} \Bigg) y - m y $$
	
	and 
	
	\begin{eqnarray*}
		\frac{\partial}{\partial x} f(x,y) &=& 1 - \frac{2 x}{\gamma} - \frac{2 x y (1 + \alpha \xi + \epsilon y)}{(1 + x^2 + \alpha \xi + \epsilon y)^2} \\
		\frac{\partial}{\partial y} f(x,y) &=& -\frac{x^2 (1+ x^2+\alpha \xi)}{(1 + x^2 + \alpha \xi + \epsilon y)^2} \\
		\frac{\partial}{\partial x} g(x,y) &=& \frac{2 \delta x y (1+(\alpha - 1) \xi + \epsilon y)}{(1+x^2 + \alpha \xi +  \epsilon y)^2} \\
		\frac{\partial}{\partial y} g(x,y) &=& \frac{\delta (x^2 + \xi) (1+x^2 + \alpha \xi)}{(1+x^2+\alpha \xi + \epsilon y)^2} - m
	\end{eqnarray*}
	
	At the trivial equilibrium point $E_0 = (0,0)$, we obtain the jacobian as 
	
	$$ J_{(0,0)} = \left[
	\begin{matrix}
		1  & 0 \\
		0 & \frac{\delta \xi - m (1+ \alpha \xi)}{1+\alpha \xi}
	\end{matrix}
	\right]$$
	
	The eigen values of this jacobian matrix are $1, \frac{\delta \xi - m (1+ \alpha \xi)}{1+\alpha \xi}$. If $\delta \xi - m (1+ \alpha \xi) > 0$, then both the eigen values have same signs. This makes the equilibrium point $E_0 = (0,0)$ unstable. If $\delta \xi - m (1+ \alpha \xi) < 0$, then both the eigen values will have opposite signs. This makes the point $(0,0)$ a saddle point.
	
	Also note that $(0,0)$ is an equilibrium point for the initial system (\ref{midm03}). The eigen values with the associated jacobian matrix are $1,-m$. Since both the eigen values have opposie signs, the point $(0,0)$ is a saddle point. 
	
	At the axial equilibrium point $E_1 = (\gamma ,0)$, we obtain the jacobian as 
	
	$$ J_{(\gamma ,0)} = \left[
	\begin{matrix}
		-1  & \frac{-\gamma ^2}{1+ \gamma ^2 + \alpha \xi} \\
		0 & \frac{(\delta - m) \gamma ^2 + \delta \xi - m (1+\alpha \xi)}{\gamma ^2 + 1 + \alpha \xi}
	\end{matrix}
	\right]$$
	
	The eigen values of this jacobian matrix are $-1, \frac{(\delta - m) \gamma ^2 + \delta \xi - m (1+\alpha \xi)}{\gamma ^2 + 1 + \alpha \xi}$. If $(\delta - m) \gamma ^2 + \delta \xi - m (1+\alpha \xi) > 0$, then both the eigenvalues will have opposite sign resulting in a saddle point. If $(\delta - m) \gamma ^2 + \delta \xi - m (1+\alpha \xi) < 0$, then it will be an asymptotically stable node. 
	
	In the absence of additional food, $E_1 = (\gamma, 0)$ remains an equilibrium point for the system (\ref{midm03}). In that case, if $(\delta - m) \gamma ^2 > m $, $E_1$ is a saddle point. Else it is a stable equilibrium point. 
	
	We now consider another axial equilibrium point which exists only for the additional food provided system (\ref{midm3}) but not for the initial system (\ref{midm03}). 
	
	At this axial equilibrium point $E_2 = (0, \frac{\delta \xi - m (1+\alpha \xi)}{m \epsilon})$, the associated jacobian matrix is given as 
	
	$$ J_{(0, \frac{\delta \xi - m (1+\alpha \xi)}{m \epsilon})} = \left[
	\begin{matrix}
		1  & 0 \\
		0 & -\frac{m (\delta \xi - m (1+\alpha \xi)) }{\delta \epsilon} 
	\end{matrix}
	\right] $$
	
	Since this equilibrium point exists in positive $xy$-quadrant only when $\delta \xi - m (1+\alpha \xi) > 0$, the eigen values of the associated jacobian matrix are of opposite side which resulting in a saddle equilibrium point. 
	
	\subsection{Interior Equilibrium point for Initial System}
	
	At the co existing equilibrium point $\bar{E} = (\bar{x}, \bar{y})$ of the initial system (\ref{midm03}) in the absence of additional food, we obtain the jacobian as 
	
	$$ J_{(\bar{x},\bar{y})} = \left[
	\begin{matrix}
		\frac{\partial}{\partial x} f(x,y)  & \frac{\partial}{\partial y} f(x,y) \\
		\frac{\partial}{\partial x} g(x,y) & \frac{\partial}{\partial y} g(x,y)
	\end{matrix}
	\right] \bigg|_{(\bar{x},\bar{y})} $$
	
	The associated characteristic equation is given by
	
	\begin{equation}
		\lambda ^2 - \text{Tr } J \bigg|_{(\bar{x},\bar{y})} \lambda + \text{Det } J \bigg|_{(\bar{x},\bar{y})} = 0.
	\end{equation}
	
	Now 
	
	\begin{eqnarray*}
		\text{Det } J \bigg|_{(\bar{x},\bar{y})} &=& \bigg( \frac{\partial f}{\partial x} \frac{\partial g}{\partial y} - \frac{\partial f}{\partial y} \frac{\partial g}{\partial x} \bigg)  \bigg|_{(\bar{x},\bar{y})} \\
		&=& \bigg( \bigg(1 - \frac{2 x}{\gamma} - \frac{2 x y (1 + \epsilon y)}{(1 + x^2 + \epsilon y)^2}  \bigg) \bigg( \frac{\delta x^2 (1+x^2)}{(1+x^2+\epsilon y)^2} - m \bigg) \bigg) \bigg|_{(\bar{x},\bar{y})} \\
		& & \ \ \ - \bigg( \bigg( -\frac{x^2 (1+ x^2)}{(1 + x^2 + \epsilon y)^2} \bigg) \bigg( \frac{2 \delta x y (1+\epsilon y)}{(1+x^2 + \epsilon y)^2}  \bigg) \bigg) \bigg|_{(\bar{x},\bar{y})} \\
	\end{eqnarray*}
	
	Since $1 - \frac{\bar{x}}{\gamma} = \frac{\bar{x} \bar{y}}{1+\bar{x}^2+\epsilon \bar{y}}$ and $\frac{\delta \bar{x}^2}{1+\bar{x}^2 + \epsilon \bar{y}} = m$, we have 
	\begin{eqnarray*}	 
		\text{Det } J \bigg|_{(\bar{x},\bar{y})} &=& \frac{m \epsilon \bar{y} + 2 m \big(1 - \frac{\bar{x}}{\gamma}\big)}{1+\bar{x}^2 + \epsilon \bar{y} } \\
		&=& \frac{m \epsilon \bar{y} + 2 m \frac{\bar{x} \bar{y}}{1+\bar{x}^2 + \epsilon \bar{y}}}{1+\bar{x}^2 + \epsilon \bar{y}} > 0
	\end{eqnarray*}
	
	Hence we have $\text{Det } J \bigg|_{(\bar{x},\bar{y})}  > 0$. This tells that both the eigen values are of the same sign. Now, in order to prove the local stability of the interior equilibrium point $(\bar{x}, \ \bar{y})$, the trace of the eigen values should be negative. This tells that both the eigen values are negative which proves the local asymptotic stability of interior equilibrium point. We now have
	
	\begin{eqnarray*}
		\text{Tr } J \bigg|_{(\bar{x},\bar{y})}  &=& \frac{\partial f}{\partial x} \bigg|_{(\bar{x},\bar{y})} + \frac{\partial g}{\partial y} \bigg|_{(\bar{x},\bar{y})} \\
		&=&  \Bigg(1 - \frac{2 x}{\gamma} - \frac{2 x y (1 + \epsilon y)}{(1 + x^2 + \epsilon y)^2}  + \frac{\delta x^2 (1+x^2)}{(1+x^2+\epsilon y)^2} - m\Bigg)  \bigg|_{(\bar{x},\bar{y})}
	\end{eqnarray*}
	
	Since $1 - \frac{\bar{x}}{\gamma} = \frac{\bar{x} \bar{y}}{1+\bar{x}^2+\epsilon \bar{y}}$ and $\frac{\delta \bar{x}^2}{1+\bar{x}^2 + \epsilon \bar{y}} = m$, we have 
	
	\begin{eqnarray*}
		\text{Tr } J \bigg|_{(\bar{x},\bar{y})}  &=& \Bigg(1 - \frac{2 x}{\gamma} - \frac{2 \big(1-\frac{x}{\gamma} \big)  (1 + \epsilon y)}{1 + x^2 + \epsilon y}  - \frac{m \epsilon y}{1+x^2+\epsilon y} \Bigg)  \bigg|_{(\bar{x},\bar{y})}
	\end{eqnarray*}
	
	From this equation, it is observed that when $1 - \frac{2 \bar{x}}{\gamma} < 0 \implies \bar{x} > \frac{\gamma}{2}$, trace of the jacobian is always negative. 
	
	Since $\bar{y} = \frac{(\delta - m) \bar{x}^2 - m}{m \epsilon}$, we have 
	\begin{eqnarray*}
		1 + \epsilon \bar{y} &=& 1 + \epsilon \bigg( \frac{(\delta - m) \bar{x}^2 - m}{m \epsilon} \bigg) = \Bigg(\frac{\delta}{m} - 1 \Bigg) \bar{x}^2 \\
		1 + \bar{x}^2 + \epsilon \bar{y} &=& 1 + \bar{x}^2 + \epsilon \bigg( \frac{(\delta - m) \bar{x}^2 - m}{m \epsilon} \bigg) = \frac{\delta}{m} \bar{x}^2 
	\end{eqnarray*}
	
	By substituting the above equations, we have
	
	\begin{equation}
		\text{Tr } J \bigg|_{(\bar{x},\bar{y})}  = \Bigg( \frac{2 m }{\delta} - 1  \Bigg) - \frac{2 m \bar{x}}{\delta \gamma} - \frac{m \epsilon \bar{y} }{1 + \bar{x}^2 + \epsilon \bar{y}} \label{tra3}
	\end{equation} 
	
	Hence $\text{Tr } J \bigg|_{(\bar{x},\bar{y})}  < 0$ when $\Bigg( \frac{2 m }{\delta} - 1  \Bigg) < 0 \implies \delta > 2 m$. Hence we have the following lemma. 
	
	\begin{lem}
		A sufficient condition for the local asymptotical stability of interior equilibrium point is $\delta > 2 m$.
	\end{lem}
	
	We now investigate the case when $m < \delta < 2 m$.
	
	We already have 
	
	$$\frac{1 - \frac{\bar{x}}{\gamma}}{\bar{x} \bar{y}} = \frac{1}{1+\bar{x}^2+\epsilon \bar{y}}$$. 
	
	Substituting this equation in (\ref{tra3}), we have
	
	\begin{eqnarray*}
		\text{Tr } J \bigg|_{(\bar{x},\bar{y})}  &=& \Bigg( \frac{2 m }{\delta} - 1  \Bigg) - \frac{2 m \bar{x}}{\delta \gamma} - \frac{m \epsilon \bar{y} }{1 + \bar{x}^2 + \epsilon \bar{y}} \\
		&=& \Bigg( \frac{2 m }{\delta} - 1  \Bigg) - \frac{2 m \bar{x}}{\delta \gamma} - m \epsilon \bar{y} \bigg( \frac{1 - \frac{\bar{x} }{\gamma}}{\bar{x} \bar{y}} \bigg) \\
		&=& \frac{-1}{\delta \gamma \bar{x}} \bigg( 2 m \bar{x}^2 - \big( (2 m - \delta) \gamma + m \epsilon \delta \big) \bar{x} + m \epsilon \delta \gamma \bigg)
	\end{eqnarray*}
	
	Let $a = 2 m, \  \ b = - \big( (2 m - \delta) \gamma + m \epsilon \delta \big), \ \ c = m \epsilon \delta \gamma$. Since $m < \delta < 2 m$, we have $ a >0,\  b <0, \ c>0$. The roots of this cubic equation $a \bar{x}^2 + b \bar{x} + c = 0 $ are
	
	$$ \bar{x} = \frac{- b \pm \sqrt{b^2 - 4 a c}}{2 a}$$
	
	If discriminant $= b^2 - 4 a c < 0$, the equation will have complex roots and the real part of the trace will be negative. This gives the asymptotic stability of the interior equilibrium point. 
	
	If discriminant $= b^2 - 4 a c > 0$, the equation will have two real roots. Since the sum of these two roots $\bar{x}_1 + \bar{x}_2 = \frac{-b}{a} > 0$ and the product of the roots $\bar{x}_1 \bar{x}_2 = \frac{c}{a} > 0$, the two real roots are positive. Since the trace at $\bar{x} = 0$ is negative, we have the positive trace from $\bar{x}_1 < \bar{x} < \bar{x}_2$ and negative otherwise. Based on the above discussion, We now have the following theorem describing the stability of interior equilibrum point. 
	
	\begin{thm}
		\begin{enumerate} The interior equilibrium point of the initial system (\ref{midm03}) 
			\item[a. ] asymptotically stable when $\delta > 2 m$.
			\item[b. ] When $m < \delta < 2 m$, 
			\begin{itemize}
				\item[i. ] If $b^2 - 4 a c > 0$ and $\bar{x}_1 < \bar{x} < \bar{x}_2$, then the interior equilibrium point is unstable.
				\item[ii. ] Else, it will be asymptotically stable.  
			\end{itemize} 
		\end{enumerate}
	\end{thm}
	
	\subsection{Interior Equilibrium point for Additional Food Provided System}
	
	At the co existing equilibrium point $E^* = (x^*, y^*)$ of the additional food provided system (\ref{midm3}), we obtain the jacobian as 
	
	$$ J_{(x^*,y^*)} = \left[
	\begin{matrix}
		\frac{\partial}{\partial x} f(x,y)  & \frac{\partial}{\partial y} f(x,y) \\
		\frac{\partial}{\partial x} g(x,y) & \frac{\partial}{\partial y} g(x,y)
	\end{matrix}
	\right] \bigg|_{(x^*,y^*)} $$
	
	The associated characteristic equation is given by
	
	\begin{equation}
		\lambda ^2 - \text{Tr } J \bigg|_{(\bar{x},\bar{y})} + \text{Det } J \bigg|_{(x^*,y^*)} = 0.
	\end{equation}
	
	Now 
	
	\begin{eqnarray*}
		\text{Det } J \bigg|_{(x^*,y^*)} &=& \bigg( \frac{\partial f}{\partial x} \frac{\partial g}{\partial y} - \frac{\partial f}{\partial y} \frac{\partial g}{\partial x} \bigg)  \bigg|_{(\bar{x},\bar{y})} \\
		&=& \bigg( \bigg(1 - \frac{2 x}{\gamma} - \frac{2 x y (1 + \alpha \xi + \epsilon y)}{(1 + x^2 + \alpha \xi + \epsilon y)^2}  \bigg) \bigg( \frac{\delta (x^2 + \xi) (1+x^2+\alpha \xi)}{(1+x^2+\alpha \xi + \epsilon y)^2} - m \bigg) \bigg) \bigg|_{(x^*,y^*)} \\
		& & \ \ \ - \bigg( \bigg( -\frac{x^2 (1+ x^2 + \alpha \xi)}{(1 + x^2 + \alpha \xi + \epsilon y)^2} \bigg) \bigg( \frac{2 \delta x y (1+(\alpha -1 ) \xi + \epsilon y)}{(1+x^2 + \alpha \xi + \epsilon y)^2}  \bigg) \bigg) \bigg|_{(x^*,y^*)} \\
	\end{eqnarray*}
	
	Since $1 - \frac{x^*}{\gamma} = \frac{x^* y^*}{1+(x^*)^2+\alpha \xi + \epsilon y^*}$ and $\frac{\delta ((x^*)^2 + \xi)}{1+(x^*)^2 + \alpha \xi + \epsilon y^*} = m$, we have 
	\begin{eqnarray*}	 
		\text{Det } J \bigg|_{(x^*,y^*)} &=& \bigg( \bigg(1 - \frac{2 x}{\gamma} - \frac{2 \bigg(1 - \frac{x}{\gamma}\bigg) (1 + \alpha \xi + \epsilon y)}{1 + x^2 + \alpha \xi + \epsilon y}  \bigg) \bigg( \frac{m (1+x^2+\alpha \xi)}{1+x^2+\alpha \xi + \epsilon y} - m \bigg) \bigg) \bigg|_{(x^*,y^*)} \\
		& & \ \ \ + \bigg( \bigg( \frac{2 \bigg( m - \frac{\delta \xi}{1 + x^2 + \alpha \xi + \epsilon y} \bigg) (1+ x^2 + \alpha \xi) (1+(\alpha -1 ) \xi + \epsilon y) \bigg( 1 - \frac{x}{\gamma} \bigg) }{(1 + x^2 + \alpha \xi + \epsilon y)^2} \bigg)  \bigg|_{(x^*,y^*)} \\
	\end{eqnarray*}
	
	Since $y^* = \frac{(\delta - m) (x^*)^2 + (\delta \xi - m (1 + \alpha \xi))}{m \epsilon}$, by substituting the following equations in the above equation, 
	
	\begin{eqnarray*}
		1 + \alpha \xi + \epsilon y^* &=& \bigg( \frac{\delta}{m} - 1 \bigg) (x^*)^2 + \frac{\delta \xi}{m}\\
		1 + (x^*)^2 + \alpha \xi + \epsilon y^* &=& \frac{\delta}{m}((x^*)^2 + \xi) \\
		m - \frac{\delta \xi}{1 + (x^*)^2 + \alpha \xi + \epsilon y^*} &=& m - \frac{\delta \xi}{\frac{\delta}{m}((x^*)^2 + \xi)} = \frac{m (x^*)^2}{(x^*)^2 + \xi}
	\end{eqnarray*}
	
	we get the determinant value as 
	
	\begin{eqnarray*}	 
		\text{Det } J \bigg|_{(x^*,y^*)} &=& \bigg( -m \bigg(1 - \frac{2 x}{\gamma} \bigg) + \frac{2 m \big(1 - \frac{x}{\gamma}\big) (1 + \alpha \xi + \epsilon y)}{1 + x^2 + \alpha \xi + \epsilon y}  + \frac{m (1+x^2+\alpha \xi) \big(1 - \frac{2 x}{\gamma}\big)}{1+x^2+\alpha \xi + \epsilon y} \bigg) \bigg|_{(x^*,y^*)} \\
		& & \ \ \ + \bigg( \frac{2 m \big( 1 - \frac{x}{\gamma} \big) (1+ x^2 + \alpha \xi) \bigg[ (1 + x^2 + \alpha \xi + \epsilon y) \big(\frac{x^2}{x^2 + \xi} - 1 \big) - \frac{\xi x^2}{x^2 + \xi} \bigg]}{(1 + x^2 + \alpha \xi + \epsilon y)^2} \bigg)  \bigg|_{(x^*,y^*)} \\
		&=& \bigg( m \big( 1 - \frac{2x}{\gamma} \big) \bigg( \frac{- \epsilon y}{1 + x^2 + \alpha \xi + \epsilon y} + \frac{2 m \big( 1 - \frac{x}{\gamma} \big) (1 +\alpha \xi + \epsilon y)}{1 + x^2 + \alpha \xi + \epsilon y} \bigg) \bigg|_{(x^*,y^*)} \\
		& & - \bigg( \frac{2 m \xi \big(1 - \frac{x}{\gamma} \big) \bigg( \frac{1 + x^2 + \alpha \xi}{x^2 + \xi}}{1 + x^2 + \alpha \xi + \epsilon y} \bigg) \bigg|_{(x^*,y^*)} \\
		&=& \frac{m}{1 + x^2 + \alpha \xi + \epsilon y} \bigg[ \epsilon y + 2 \bigg( 1 - \frac{x}{\gamma} \bigg) \bigg( \frac{(1 + (\alpha - 1) \xi) x^2}{x^2 + \xi} \bigg) \bigg] \bigg|_{(x^*,y^*)}  \\
		&=& \frac{\epsilon \delta y^* ((x^*)^2 + \xi) + 2 \delta \big(1 - \frac{x^*}{\gamma}\big) (1 + (\alpha - 1) \xi) (x^*)^2}{(1 + (x^*)^2 + \alpha \xi + \epsilon y^*)^2} 
	\end{eqnarray*}
	
	In the absence of mutual interference ($\epsilon = 0$), the determinant is positive when $1 + (\alpha -1) \xi > 0$. This is in line with the results obtained in \cite{V3JTB}. 
	
	In the presence of mutual inerference and the absence of additional food, the determinant is always positive. This is in line with the results obtained in the previous subsection of this paper. 
	
	In the presence of mutual interference and additional food, the determinant value is not always positive. The determinant is positive when $1 + (\alpha -1) \xi > 0$. The term $1 + (\alpha -1) \xi$ is positive for all $\alpha > 1 $ and $\xi > 1$. 
	
	\begin{figure}[ht]
		\centering
		\includegraphics[width=0.8\textwidth]{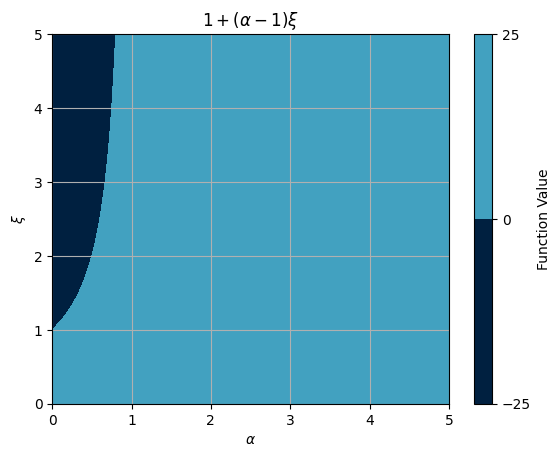}
		\caption{The contour plot depicting the value of $1 + (\alpha - 1) \xi$ for a combination of positive-$\alpha \xi$. }
		\label{expr3}
	\end{figure}
	
	The figure \ref{expr3} depicts the region where $1 + (\alpha -1) \xi > 0$. The light blue region gives $1 + (\alpha - 1) \xi > 0$ resulting in eigen values having same sign. The dark blue region where $\alpha < 1$ and $\xi > 1$ gives negative value for the expression $1 + ( \alpha - 1) \xi$. In this region, the derivative can still remain positive if 
	\begin{equation} \label{epsdel}
		\epsilon > \frac{- 2 }{y} \bigg( 1 - \frac{x}{\gamma} \bigg) ( 1 + (\alpha - 1) \xi ) > 0
	\end{equation}
	
	Derivative will be negative for the remaining values of $\epsilon$. Hence the equilibrium point will be saddle for these remaining values of $\epsilon$.
	
	The equation (\ref{epsdel}) can be simplified as a quartic equation, 
	$$ \bigg(\frac{\delta}{m} - 1 \bigg) (x^*)^4 - \frac{2}{\gamma} \bigg(1 + (\alpha-1) \xi \bigg) (x^*)^3 + (1 + (\alpha - 1) \xi + 2 \xi \bigg(\frac{\delta}{m} - 1 \bigg)) (x^*)^2 + \xi (\frac{\delta \xi}{m} - (1 + \alpha \xi)) > 0$$
	
	From the theory of equations, this quartic expression has a negative sum of its roots and positive product of its roots. Therefore, it can have atmost two positive roots $x_1,\ x_2$. Also, this quartic expression is positive when $x^* = 0, \frac{\gamma}{2}, \ \gamma$. Therefore, if $0 < x_1 < x_2 < \gamma$, the derivative will be negative in the region $x_1 < x^* < x_2$ resulting in a saddle point. Else, it will be positive. In which case, a further analysis on the trace of jacobian will reveal whether $(x^*,y^*)$ is a stable or unstable equilibrium point.

	The trace of the jacobian matrix is given by 
	
	\begin{eqnarray*}
		\text{Tr } J \bigg|_{(x^*,y^*)}  &=& \frac{\partial f}{\partial x} \bigg|_{(x^*,y^*)} + \frac{\partial g}{\partial y} \bigg|_{(x^*,y^*)} \\
		&=&  \Bigg(1 - \frac{2 x}{\gamma} - \frac{2 x y (1 + \alpha \xi +  \epsilon y)}{(1 + x^2 + \alpha \xi + \epsilon y)^2}  + \frac{\delta (x^2 + \xi) (1+x^2+\alpha \xi)}{(1+x^2+\alpha \xi + \epsilon y)^2} - m\Bigg)  \bigg|_{(x^*,y^*)}
	\end{eqnarray*}
	
	Since $1 - \frac{x^*}{\gamma} = \frac{x^* y^*}{1+(x^*)^2+\alpha \xi + \epsilon y^*}$ and $\frac{\delta ( (x^*)^2 + \xi)}{1+(x^*)^2 + \alpha \xi + \epsilon y^*} = m$, we have 
	
	\begin{eqnarray*}
		\text{Tr } J \bigg|_{(x^*,y^*)}  &=&  \Bigg(1 - \frac{2 x}{\gamma} - \frac{2 \bigg(1 - \frac{x}{\gamma}\bigg) (1 + \alpha \xi +  \epsilon y)}{1 + x^2 + \alpha \xi + \epsilon y}  + \frac{m (1+x^2+\alpha \xi)}{1+x^2+\alpha \xi + \epsilon y} - m\Bigg)  \bigg|_{(x^*,y^*)} \\
		&=&  \Bigg( \frac{\big(1-\frac{2 x}{\gamma} \big) x^2 -  (1 + \alpha \xi + \epsilon y)}{1 + x^2 + \alpha \xi + \epsilon y}  - \frac{m \epsilon y}{1+x^2+\alpha \xi + \epsilon y} \Bigg)  \bigg|_{(x^*,y^*)}
	\end{eqnarray*}
	
	From the above equation, it is observed that the trace of the jacobian will be negative when $x^*=0$ and $1-\frac{2 x^*}{\gamma} < 0 \implies x^* > \frac{\gamma}{2}$. When $0 < x^* < \frac{\gamma}{2}$, the trace of the jacbian can be calculated by substituting $y^* = \frac{(\delta - m) {x^*}^2 + \delta \xi - m (1+\alpha \xi)}{m \epsilon}$ in the above equation. Upon simplification, we have 
	\begin{equation}
		\text{Tr } J \bigg|_{(x^*,y^*)}  = \frac{- \frac{2 {x^*}^3}{\gamma} +  \bigg( 2 - \frac{\delta}{m} - (\delta - m) \bigg) {x^*}^2 - \frac{\delta \xi}{m} - (\delta \xi - m (1+\alpha \xi))}{1 + x^* + \alpha \xi + \epsilon y^*} \label{tra4}
	\end{equation} 
	
	Since $y^* > 0$, we have $(\delta - m) {x^*}^2  - (\delta \xi - m (1+\alpha \xi) > 0$. 
	
	\begin{eqnarray*}
		\text{Tr } J \bigg|_{(x^*,y^*)}  &=& \frac{- \frac{2 {x^*}^3}{\gamma} +  \bigg( 2 - \frac{\delta}{m} - (\delta - m) \bigg) {x^*}^2 - \frac{\delta \xi}{m} - (\delta \xi - m (1+\alpha \xi))}{1 + x^* + \alpha \xi + \epsilon y^*} \\
		&<& \frac{- \frac{2 {x^*}^3}{\gamma} +  \bigg( 2 - \frac{\delta}{m}  \bigg) {x^*}^2 - \frac{\delta \xi}{m} }{1 + x^* + \alpha \xi + \epsilon y^*} 
	\end{eqnarray*} 
	Hence $\text{Tr } J \bigg|_{(x^*,y^*)}  < 0$ when $\Bigg( 2 - \frac{\delta}{m}  \Bigg) < 0 \implies \delta > 2 m$. Hence we have the following lemma. 
	
	\begin{lem}
		A sufficient condition for the local asymptotical stability of interior equilibrium point is $\delta > 2 m$.
	\end{lem}
	
	Now, we are left to find the sign of trace of jacobian when $m < \delta < 2 m$ and $0 < x^* < \frac{\gamma}{2}$.
	
	The value of expression $- \frac{2 {x^*}^3}{\gamma} +  \bigg( 2 - \frac{\delta}{m}  \bigg) {x^*}^2 - \frac{\delta \xi}{m} $ at $x^*=0$ and $x^*=\frac{\gamma}{2}$ are negative. Hence it can have either two or no positive roots in this interval. From sturm's sequence for cubic equations, this expression exhibits two positive roots between $0$ and $\frac{\gamma}{2}$ only when $$\frac{18 \delta \xi }{m \gamma^2} < \bigg(2 - \frac{\delta}{m} \bigg)^2.$$ 
	
	Now let us consider the case $\delta \xi - m(1+\alpha \xi) < 0$. From the study of isoclines in the previous section, we derived in (\ref{condiso3}) that the interior equilibrium point exists when $\delta \xi - m(1+\alpha \xi) > - (\delta - m) \gamma^2$. From (\ref{tra4}), we have
	
	\begin{eqnarray*}
		\text{Tr } J \bigg|_{(x^*,y^*)}  &=& \frac{- \frac{2 {x^*}^3}{\gamma} +  \bigg( 2 - \frac{\delta}{m} - (\delta - m) \bigg) {x^*}^2 - \frac{\delta \xi}{m} - (\delta \xi - m (1+\alpha \xi))}{1 + x^* + \alpha \xi + \epsilon y^*} \\
		&<& \frac{- \frac{2 {x^*}^3}{\gamma} +  \bigg( 2 - \frac{\delta}{m}  \bigg) {x^*}^2 - (\delta \xi - m (1+\alpha \xi)) }{1 + x^* + \alpha \xi + \epsilon y^*} 
	\end{eqnarray*} 
	
	As the product of roots is positive and the sum of product of two roots is zero in the above expression, trace of the jacobian is negative. 
	
	Now consider the case when $\delta \xi - m(1+\alpha \xi) > 0$.
	The trace is negative when $2 - \frac{\delta}{m} - (\delta - m) < 0 \implies \delta > \frac{2+m}{1+\frac{1}{m}}$. Else, trace is positive only for the $x^*$ between the two positive roots $0 < x_1 < x^* < x_2 < \frac{\gamma}{2}$ and remains negative otherwise.  
	
	From (\ref{tra4}), we have
	
	\begin{equation*}
		\text{Tr } J \bigg|_{(x^*,y^*)}  = \frac{- \frac{2 {x^*}^3}{\gamma} +  \bigg( 2 - \frac{\delta}{m} - (\delta - m) \bigg) {x^*}^2 - \frac{\delta \xi}{m} - (\delta \xi - m (1+\alpha \xi)}{1 + x^* + \alpha \xi + \epsilon y^*} 
	\end{equation*} 
	
	Let $x_1$ and $x_2$ be two positive roots of this cubic equation $a {x^*}^3 + b {x^2} + d = 0 $, we have $a = - \frac{2}{\gamma}, \  \ b =  2 - \frac{\delta}{m} - (\delta - m), \ \ c = 0, \ \ d = - \frac{\delta \xi}{m} - (\delta \xi - m (1+\alpha \xi)$.
	
	Let $\Omega = \{(\alpha, \xi) | \ 1 + (\alpha - 1) \xi > 0 \text{ or } 1 + (\alpha - 1) \xi < 0 \text{ and } 0 < \frac{-2}{y^*} (1 - \frac{x^*}{\gamma}) (1+(\alpha - 1) \xi) < \epsilon \}$.
	
	Based on the above discussions, we have the following theorem dealing with the stability of the interior equilibrum point. 
	
	\begin{thm}
		The interior equilibrium point of the additional food provided system (\ref{midm3}) is
		\begin{enumerate} 
			\item[a. ] asymptotically stable when $(\alpha, \xi ) \in \Omega$ and $\delta > 2 m$.
			\item[b. ] When $(\alpha, \xi ) \in \Omega$ and $m < \delta < 2 m$, 
			\begin{itemize}
				\item[i. ] If $\delta \xi - m(1+\alpha \xi) > 0$, $\delta < \frac{2 + m}{1+\frac{1}{m}}$ and $0 < x_1 < x^* < x_2 < \frac{\gamma}{2}$, then the interior equilibrium point is a unstable point.
				\item[ii. ] Else, it will be asymptotically stable.  
			\end{itemize} 
			\item[c.] Saddle when $1 + (\alpha - 1) \xi < 0$ and $0 < \epsilon < \frac{-2}{y^*} (1 - \frac{x}{\gamma}) (1+(\alpha - 1) \xi) $. 
		\end{enumerate}
	\end{thm} 
	
	\section{Global dynamics of Initial System}\label{secglob1}
	
	In similar lines to \cite{V3JTB}, the initial system (\ref{midm03}) exhibits Hopf bifurcation when trace of the jacobian at interior equilbrium point is $0$,
	
	\begin{eqnarray*}
		\textit{i.e.,} \text{ Tr } J \bigg|_{(x^*,y^*)}  = & & 2 m \bar{x}^2 - \big( (2 m - \delta) \gamma + m \epsilon \delta \big) \bar{x} + m \epsilon \delta \gamma = 0, \\
		\text{ where }  & & \bar{x} = \frac{\epsilon \gamma \delta + \sqrt{(\epsilon \gamma \delta)^2 + 4 m \gamma ((\delta - m) \gamma + \delta \epsilon)}}{2 ((\delta - m)\gamma + \delta \epsilon)} \  > 0
	\end{eqnarray*}
	
	Based on the discussions so far, the parameter space of the initial system can be divided into four regions. 
	
	Let $\mathfrak{R} = \{ \gamma \ \  | \ \ 2 m \bar{x}^2 - \big( (2 m - \delta) \gamma + m \epsilon \delta \big) \bar{x} + m \epsilon \delta \gamma = 0\}$.
	
	\begin{eqnarray}
		\text{Region I} = \bigg\{ (\gamma, \epsilon, \delta, m) &|& \gamma \leq \frac{\epsilon \gamma \delta + \sqrt{(\epsilon \gamma \delta)^2 + 4 m \gamma ((\delta - m) \gamma + \delta \epsilon)}}{2 ((\delta - m)\gamma + \delta \epsilon)}  \bigg\} \label{r13} \\
		\text{Region II} = \bigg\{ (\gamma, \epsilon, \delta, m) &|& \gamma > \frac{\epsilon \gamma \delta + \sqrt{(\epsilon \gamma \delta)^2 + 4 m \gamma ((\delta - m) \gamma + \delta \epsilon)}}{2 ((\delta - m)\gamma + \delta \epsilon)} \text{ and } \nonumber \\ 
		& & \gamma \leq 4 \epsilon + \sqrt{16 \epsilon^2 + 27}\bigg\} \label{r23}\\
		\text{Region III} = \bigg\{ (\gamma, \epsilon, \delta, m) &|& \gamma > \frac{\epsilon \gamma \delta + \sqrt{(\epsilon \gamma \delta)^2 + 4 m \gamma ((\delta - m) \gamma + \delta \epsilon)}}{2 ((\delta - m)\gamma + \delta \epsilon)} \text{ and } \nonumber \\ & & 4 \epsilon + \sqrt{16 \epsilon^2 + 27} < \gamma < \mathfrak{R} \bigg\}  \label{r33} \\
		\text{Region IV} = \bigg\{ (\gamma, \epsilon, \delta, m) &|& \gamma > \frac{\epsilon \gamma \delta + \sqrt{(\epsilon \gamma \delta)^2 + 4 m \gamma ((\delta - m) \gamma + \delta \epsilon)}}{2 ((\delta - m)\gamma + \delta \epsilon)}  \nonumber \\ & & \text{ and } \gamma > \mathfrak{R} \bigg\} \label{r43} 
	\end{eqnarray}
	
\begin{figure}[ht]
		\centering
		\includegraphics[width=0.8\textwidth]{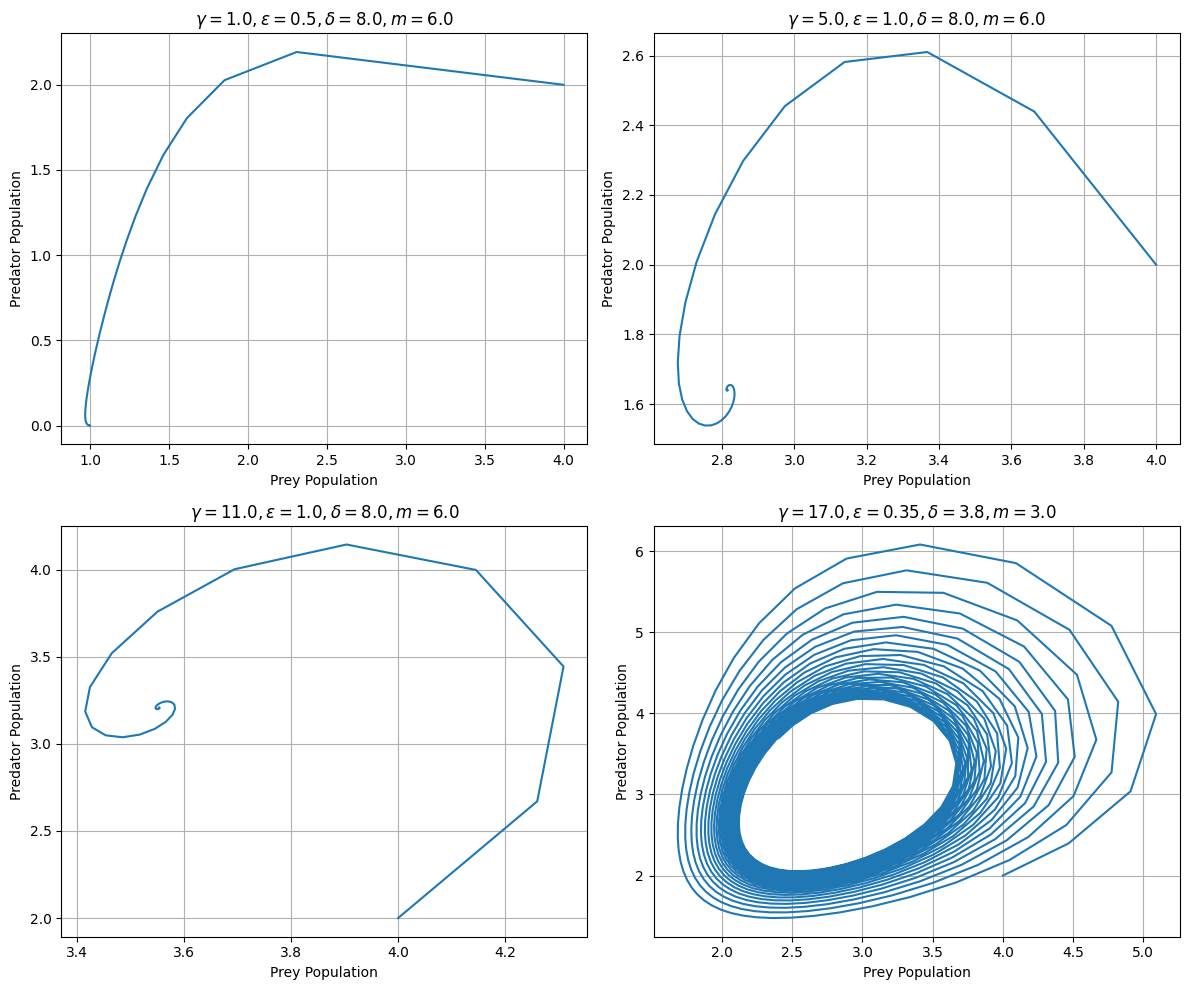}
		\caption{This figure depicts the phase portrait of the initial system (\ref{midm03}) in the regions I,II, III, IV respectively. }
		\label{init3}
	\end{figure}

	In these four regions, interior equilibrium point does not exist in region 1. It exists and is stable in regions II and III. In region IV, it is unstable and admits a limit cycle. The phase portrait in these four regions is depicted in figure \ref{init3}.
	
	\section{Global dynamics of Additional Food Provided System} \label{secglob2}
	
	In this section, we study the global dynamics of the system (\ref{midm3}) with respect to the changes in the $(\alpha, \xi)$ control parameter space. 
	
	From the stability analysis of equilibrium points, it is observed that the nature of equilibrium points $(0,0), \ (\gamma, 0),\ (x^*,y^*)$ depends on the signs of the expressions $\delta \xi -m (1+\alpha \xi)$, $(\delta - m) \gamma^2 + \delta \xi - m (1+\alpha \xi)$, $- \frac{2 {x^*}^3}{\gamma} +  \bigg( 2 - \frac{\delta}{m} - (\delta - m) \bigg) {x^*}^2 - \frac{\delta \xi}{m} - (\delta \xi - m (1+\alpha \xi))$. Based on the results obtained so far, the behavior of equilibrium points of the  system (\ref{midm3}) can be summarized in the table \ref{sumaddfood}.
	
	\begin{table}[ht]
		\centering
		\begin{tabular}{|c|c|c|c|c|}
			\hline
			& \textbf{(0,0)} & \textbf{($\gamma$,0)} & \textbf{(0,$\frac{\delta \xi - m (1+\alpha \xi)}{m \epsilon}$)} & \textbf{$(x^*,y^*)$} \\
			\hline
			$\delta \xi - m(1 + \alpha \xi)>0,\ \delta > 2m$ & Unstable & Saddle & Saddle & Stable \\
			\hline
			$0 > \delta \xi - m(1 + \alpha \xi) > -(\delta - m) \gamma^2$ & Saddle & Saddle & Not Exist & Stable \\
			\hline
			$\delta \xi - m(1 + \alpha \xi) < -(\delta - m) \gamma^2$ & Saddle & Stable & Not Exist & Not Exist \\
			\hline
		\end{tabular}
		\caption{Stability of equilibrium points of additional food provided system (\ref{midm3}).}
		\label{sumaddfood}
	\end{table}
	
	From table \ref{sumaddfood}, we consider the following curves in the positive quadrant of $(\alpha,\xi)$ space:
	
	\begin{eqnarray}
		\delta \xi - m(1+\alpha \xi) = 0 \label{PEC} \\
		\delta \xi - m(1 + \alpha \xi) = -(\delta - m) \gamma^2 \label{TBC}
	\end{eqnarray}
	
	In the positive $(\alpha,\xi)$ quadrant, when we fix the value of $\alpha$ and increase the value of $\xi$ from (\ref{TBC}), we observe that the existence of the interior equilibrium point. This results in the pseudo transcritical bifurcation due to the exchange of stability between $(\gamma,0)$ and $(x^*,y^*)$. Here, the nature of $(0,0)$ remains unchanged. 
	
	Also, the increase of $\xi$ from (\ref{PEC}) results in a pseudo transcritical bifurcation due to the exchange of stability from $(0,0)$ and $(0,\frac{\delta \xi - m (1+\alpha \xi)}{m \epsilon})$. The curve \ref{PEC} is observed to be a prey elimination curve. 
	
	\subsection{Effect of quality of additional food \texorpdfstring{($\alpha$)}{ }}
	
	\begin{figure}[ht]
		\centering
		\includegraphics[width=\textwidth]{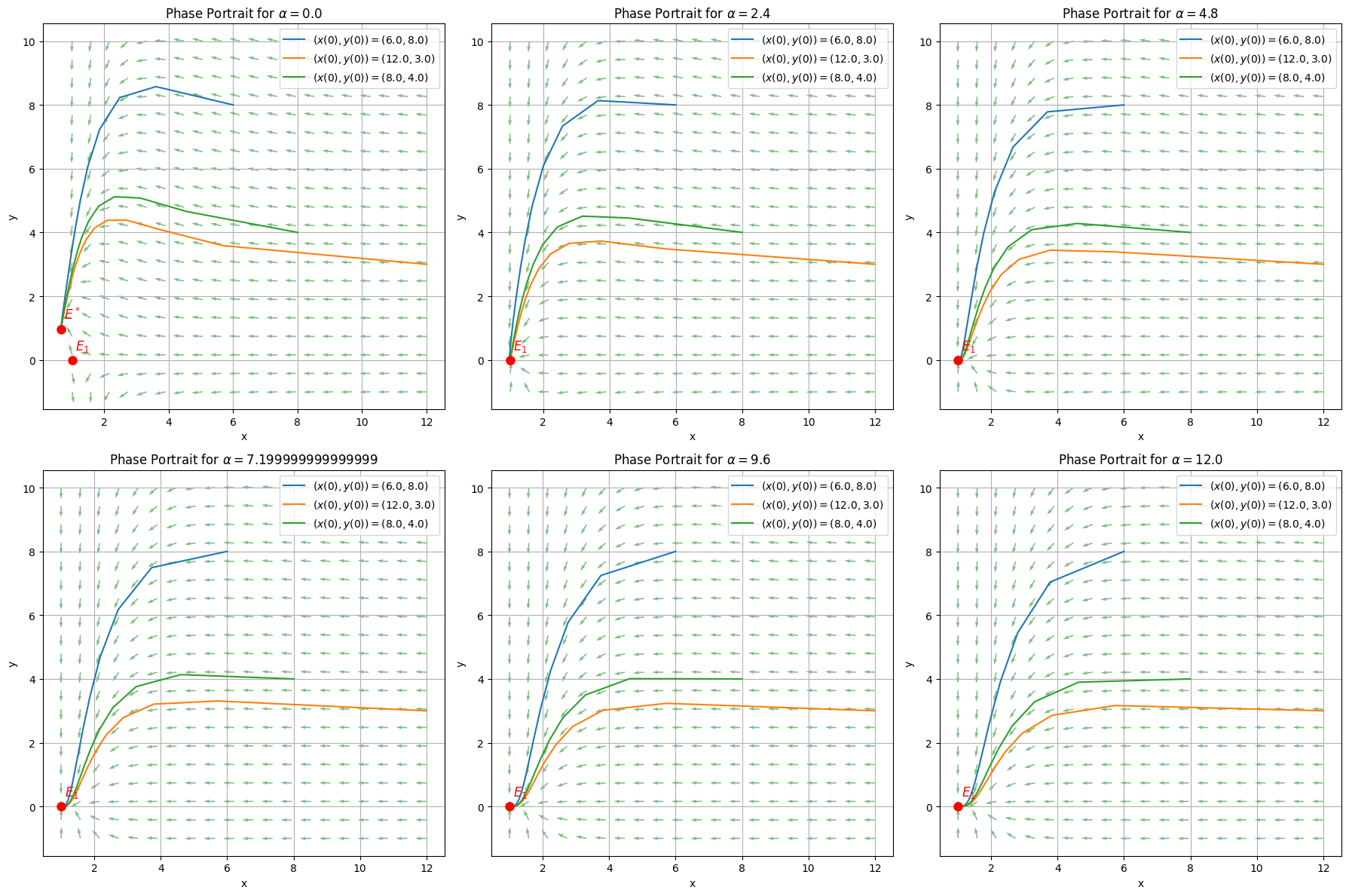}
		\caption{Phase portrait of additional food provided system (\ref{midm3}) for an array of values of $\alpha$. Here the remaining parameters belong to Region I (\ref{r13}).}
		\label{r1alpha3}
	\end{figure}
	
	From figure \ref{r1alpha3}, it is observed that the presence of additional food is shifting the dynamical system (\ref{midm3}) from equilibrium point $E^*$ to $E_1$ in region I (\ref{r13}). It is also observed that the impact of $\xi$ is minimal in region I.  The parameter values used for this analysis are $\gamma = 1.0, \ \xi = 1.0,\  \epsilon = 0.5,\  \delta = 8.0, \ m = 6.0$.
	
	\begin{figure}[ht]
		\centering
		\includegraphics[width=\textwidth]{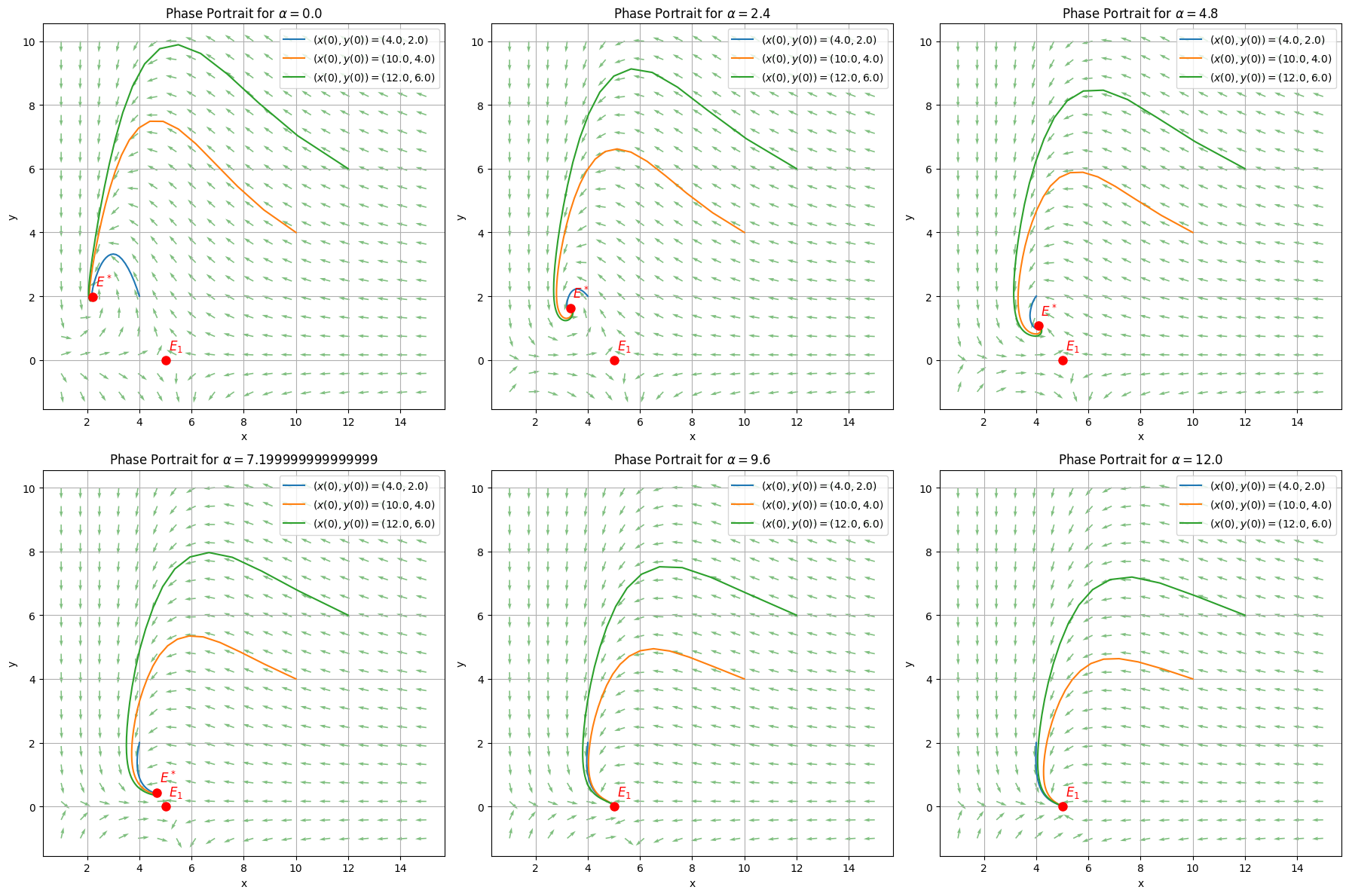}
		\caption{Phase portrait of additional food provided system (\ref{midm3}) for an array of values of $\alpha$. Here the remaining parameters belong to Region II (\ref{r23}).}
		\label{r2alpha3}
	\end{figure}
	
	From figure \ref{r2alpha3}, it is observed that the high quality of additional food is required for the dynamical system (\ref{midm3}) to shift from equilibrium point $E^*$ to $E_1$ in region II (\ref{r23}). It is also observed that the large values of $\xi$ can make the system reach $E_1$ faster in region II.  The parameter values used for this analysis are $\gamma = 5.0,\  \xi = 1.0, \ \epsilon = 1.0,\  \delta = 8.0,\  m = 6.0$.
	
	\newpage
	
	\begin{figure}[ht]
		\centering
		\includegraphics[width=\textwidth]{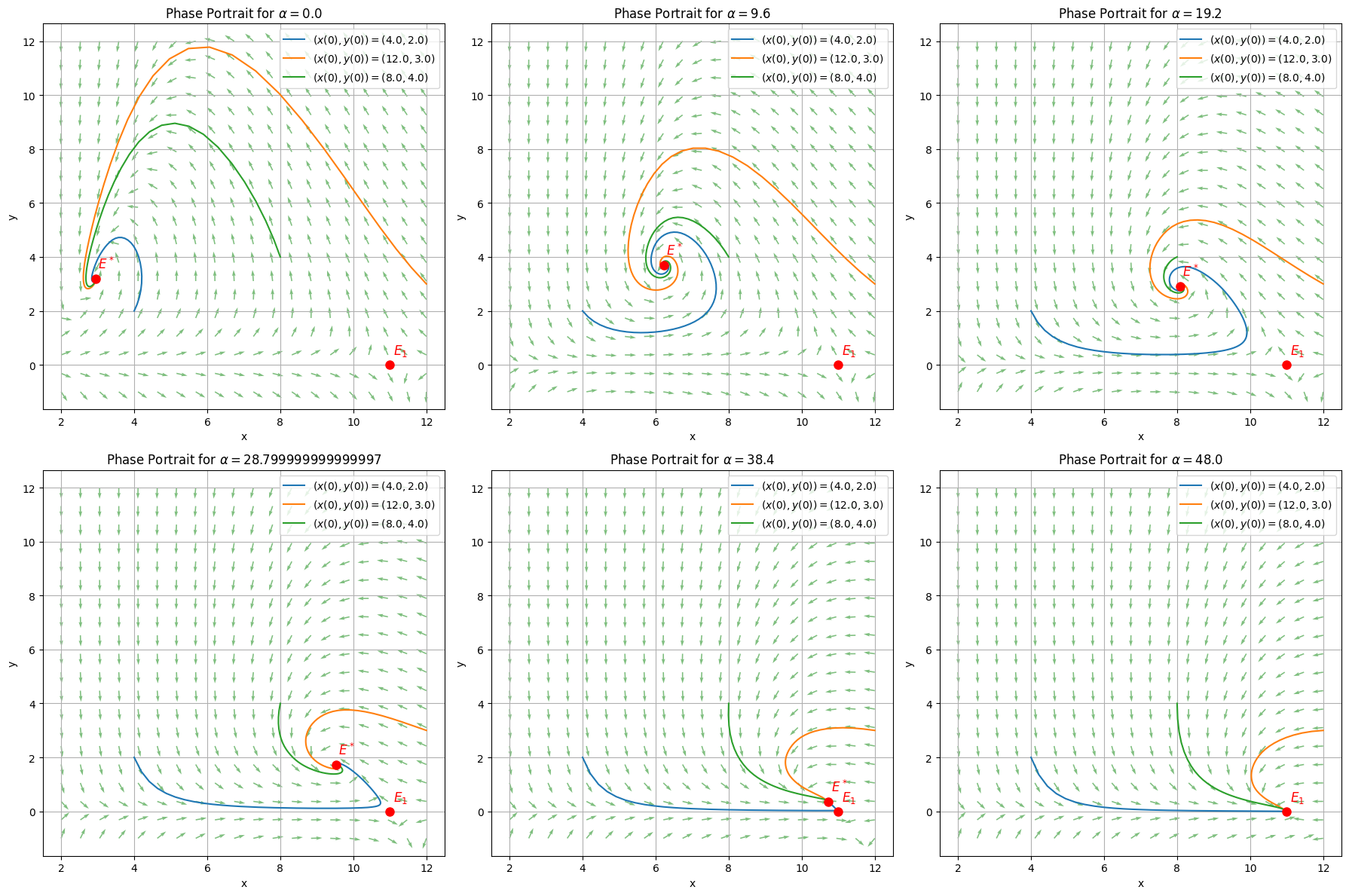}
		\caption{Phase portrait of additional food provided system (\ref{midm3}) for an array of values of $\alpha$. Here the remaining parameters belong to Region III (\ref{r33}).}
		\label{r3alpha3}
	\end{figure}
	
	From figure \ref{r3alpha3}, it is observed that the very high quality of additional food (In this case, $\alpha > 38.4$) is required for the dynamical system (\ref{midm3}) to shift from equilibrium point $E^*$ to $E_1$ in region III (\ref{r33}). It is also observed that the large values of $\xi$ can make the system reach $E_1$ faster in region III.  The parameter values used for this analysis are $\gamma = 11.0, \ \xi = 1.0, \ \epsilon = 1.0,\  \delta = 8.0, \ m = 6.0$.
	
	\newpage
	
	\begin{figure}[ht]
		\centering
		\includegraphics[width=\textwidth]{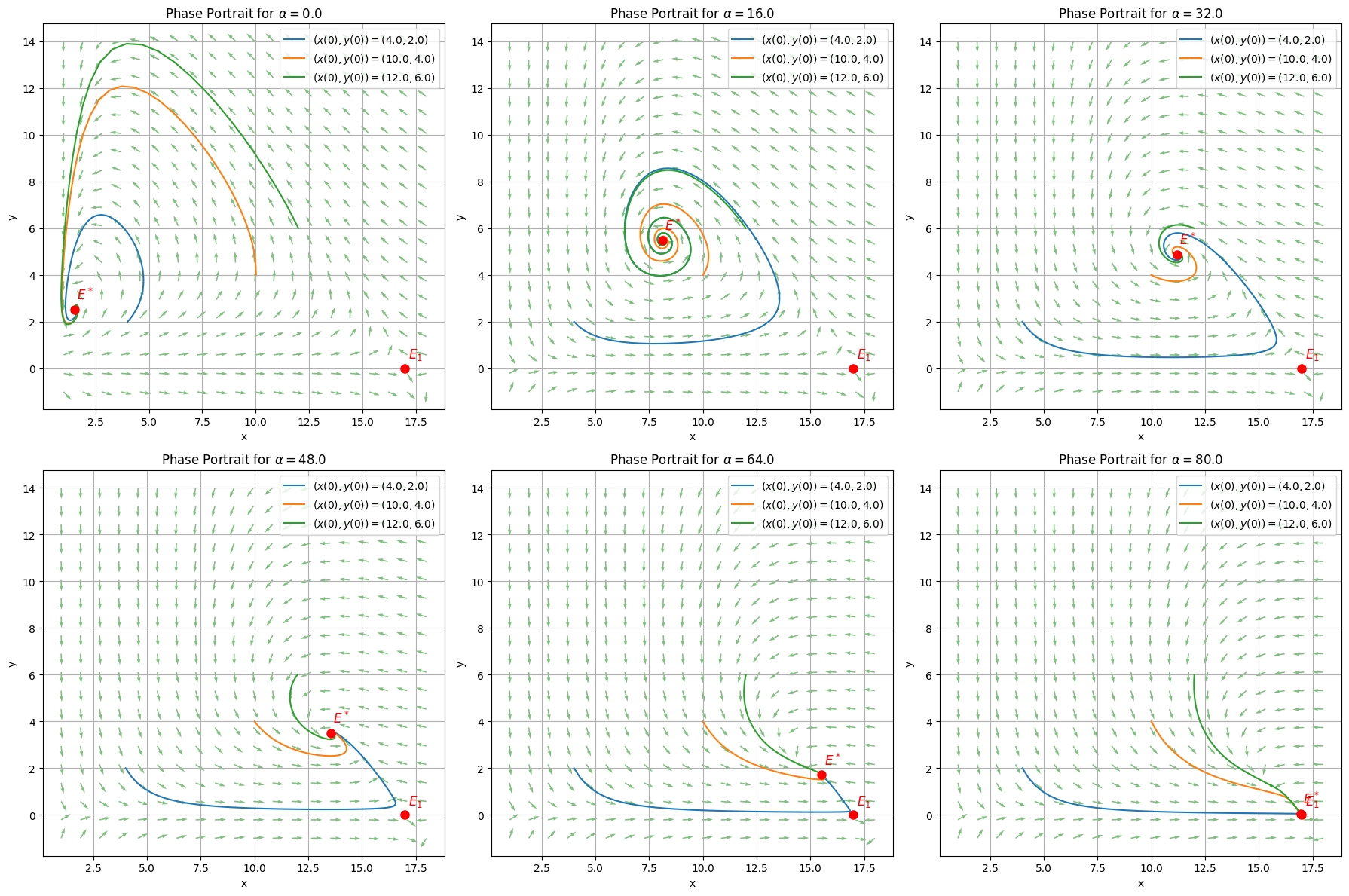}
		\caption{Phase portrait of additional food provided system (\ref{midm3}) for an array of values of $\alpha$. Here the remaining parameters belong to Region IV (\ref{r43}).}
		\label{r4alpha3}
	\end{figure}
	
	From figure \ref{r4alpha3}, it is observed that even the very high quality of additional food (In this case, $\alpha > 80.0$), the dynamical system (\ref{midm3}) is tending to the equilibrium point $E^*$ in region IV (\ref{r43}). Further increase in $\alpha$ drives the system to the equilibrium point $E_1$. It is also observed that the large values of $\xi$ can make the system reach $E_1$ faster in region IV. However, $\xi < 1$ will result in a limit cycle and the limit cycle behavior is observed even for moderately higher values of $\alpha$.  The parameter values used for this analysis are $\gamma = 17.0,\  \xi = 1.0,\  \epsilon = 0.35, \ \delta = 3.8, \ m = 3.0$.
	
	From figures \ref{r1alpha3} - \ref{r4alpha3}, the change in nature of the equilibrium points are observed with increased values of $\alpha$. There is also a possibility of bifurcation between $E_1 = (\gamma,0)$ and $E^* = (x^*,y^*)$. 
	
	\subsection{Effect of quantity of additional food \texorpdfstring{($\xi$)}{ }}
	
	\begin{figure}[ht]
		\centering
		\includegraphics[width=\textwidth]{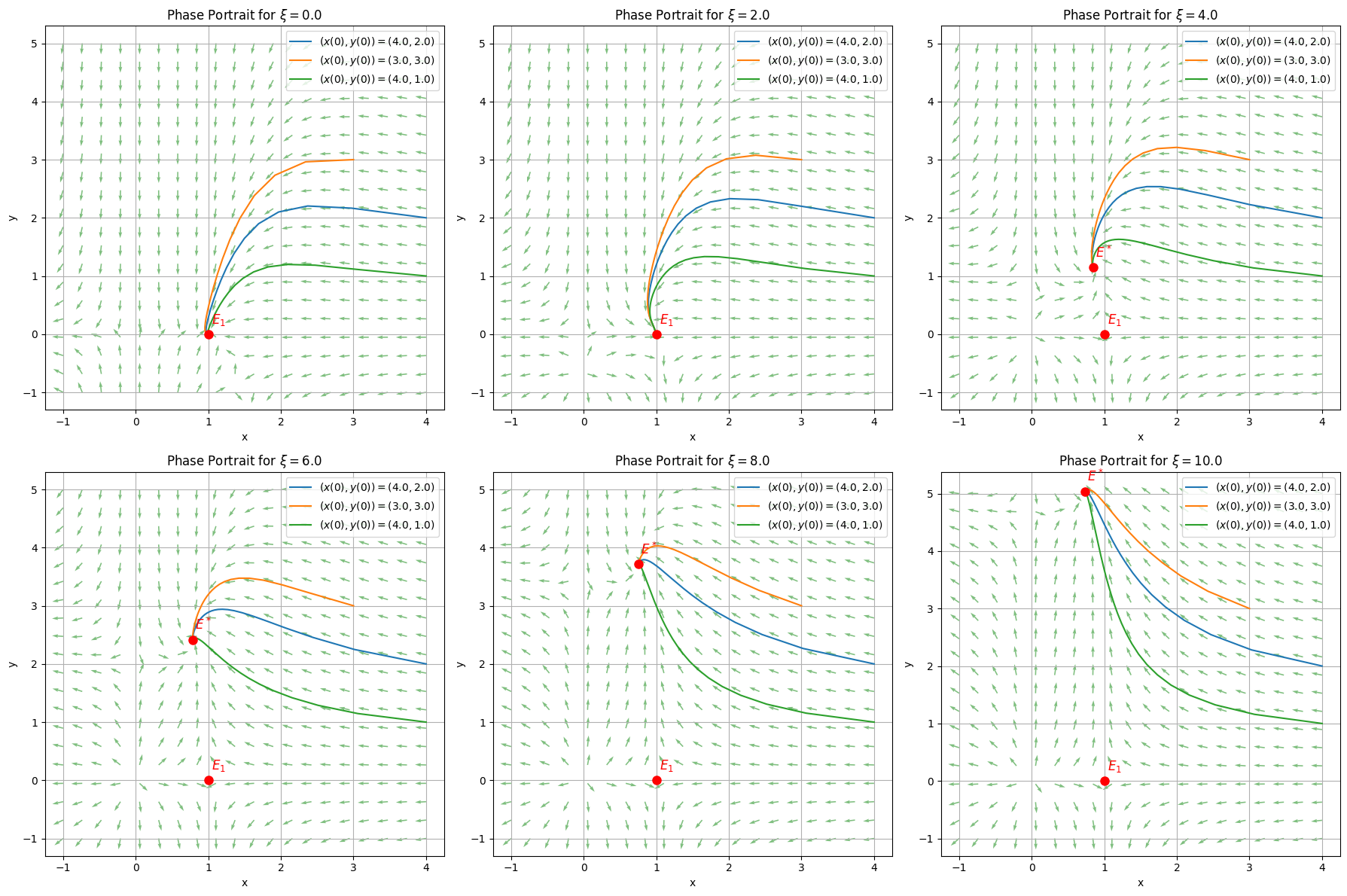}
		\caption{Phase portrait of additional food provided system (\ref{midm3}) for an array of values of $\xi$. Here the remaining parameters belong to Region I (\ref{r13}).}
		\label{r1xi3}
	\end{figure}
	
	From figure \ref{r1xi3}, it is observed that on the supply of low quantity of additional food, there is a shift in the stability of the dynamical system (\ref{midm3}) from equilibrium point $E^1$ to $E_*$ in region I (\ref{r13}). This is in contrary to the observations made in region 1 for various values of $\alpha$. It is also observed that the higher values of $\alpha$ will drive system to only $E_1$ in region I.  The parameter values used for this analysis are $\gamma = 1.0, \ \alpha = 1.0,\  \epsilon = 0.5,\  \delta = 8.0,\  m = 6.0$.
	
	\newpage
	
	\begin{figure}[ht]
		\centering
		\includegraphics[width=\textwidth]{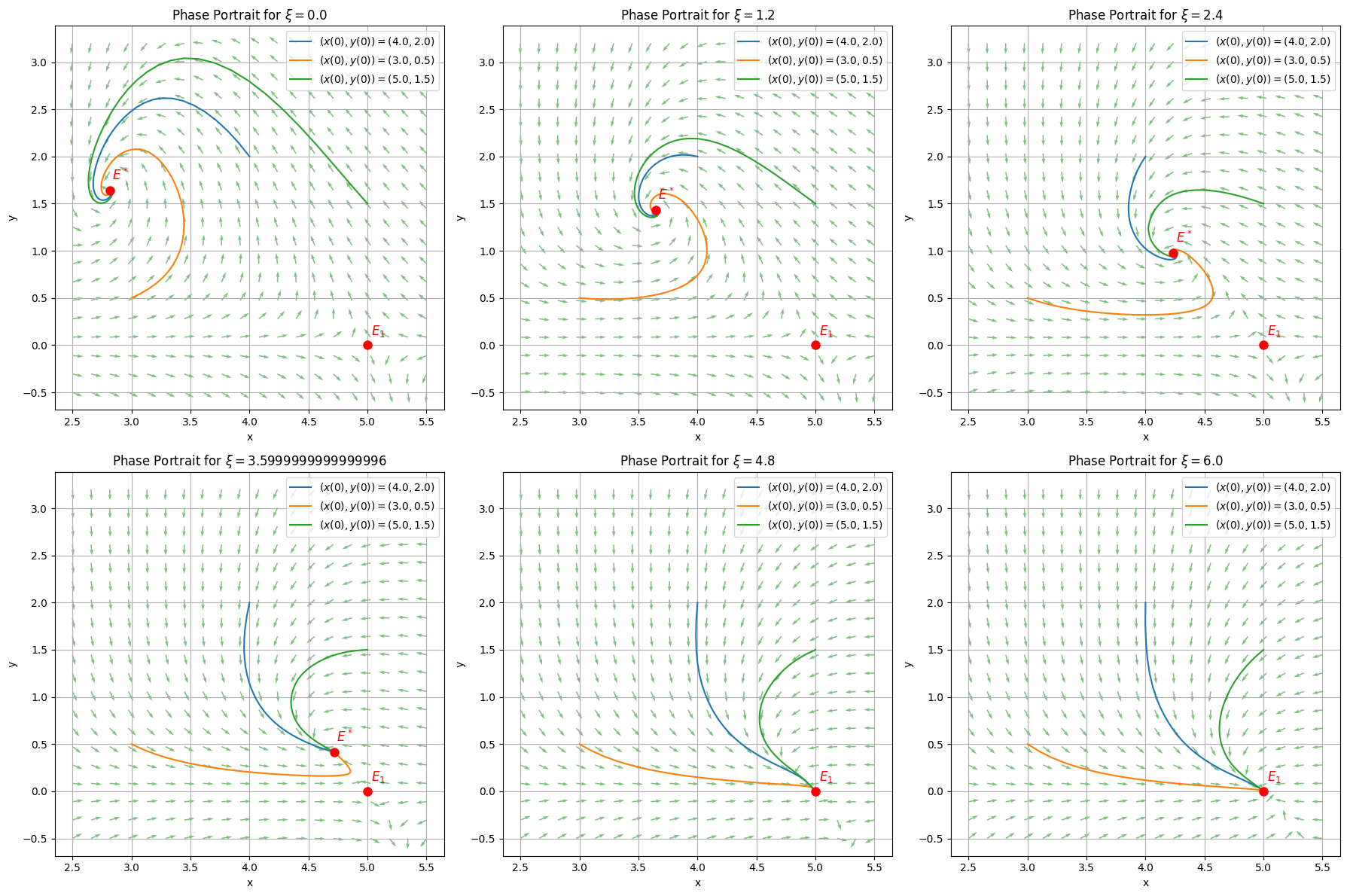}
		\caption{Phase portrait of additional food provided system (\ref{midm3}) for an array of values of $\xi$. Here the remaining parameters belong to Region II (\ref{r23}).}
		\label{r2xi3}
	\end{figure}
	
	From figure \ref{r2xi3}, it is observed that even the low quality of additional food is required for the dynamical system (\ref{midm3}) to shift from equilibrium point $E^*$ to $E_1$ in region II (\ref{r23}). It is also observed that the large values of $\alpha$ can make the system reach $E_1$ faster in region II.  The parameter values used for this analysis are $\gamma = 5.0,\  \alpha = 3.0,\  \epsilon = 1.0,\  \delta = 8.0,\  m = 6.0$.
	
	\newpage
	
	\begin{figure}[ht]
		\centering
		\includegraphics[width=\textwidth]{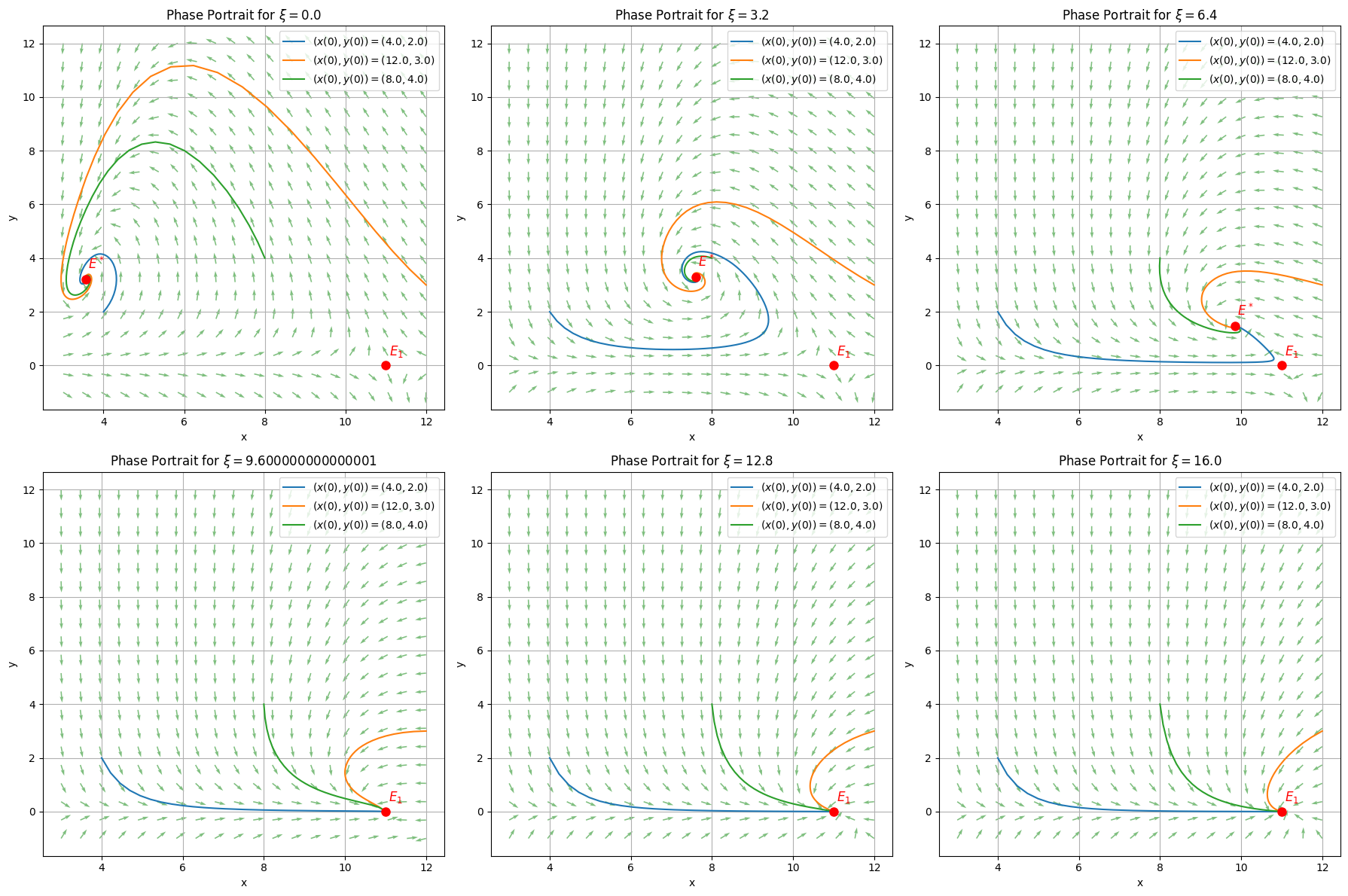}
		\caption{Phase portrait of additional food provided system (\ref{midm3}) for an array of values of $\xi$. Here the remaining parameters belong to Region III (\ref{r33}).}
		\label{r3xi3}
	\end{figure}
	
	From figure \ref{r3xi3}, it is observed that even the higher quantity of additional food will continue to drive system (\ref{midm3}) only to $E^*$. However increase in $\alpha$ can drive system faster to $E_1$ in region III. The parameter values used for this analysis are $\gamma = 11.0,\  \alpha = 6.0, \ \epsilon = 1.0,\  \delta = 8.0, \ m = 6.0$.
	
	\newpage
	
	\begin{figure}[ht]
		\centering
		\includegraphics[width=\textwidth]{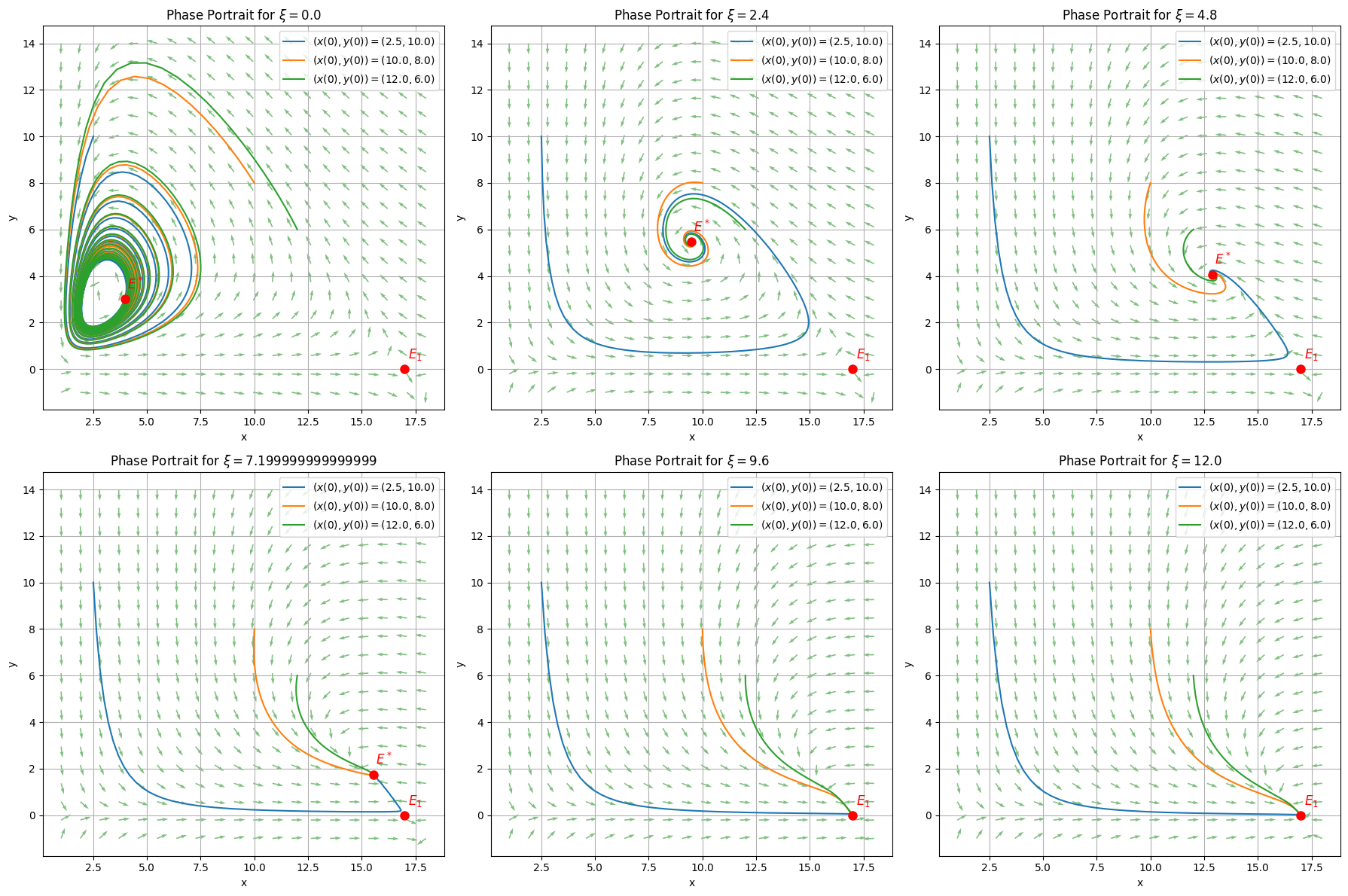}
		\caption{Phase portrait of additional food provided system (\ref{midm3}) for an array of values of $\xi$. Here the remaining parameters belong to Region IV (\ref{r43}).}
		\label{r4xi3}
	\end{figure}
	
	In the absence of additional food, the dynamical system (\ref{midm3}) exhibits an unstable limit cycle. From figure \ref{r4xi3}, it is observed that there is a shift in the stability of the system from limit cycle to $E^*$ and finally to $E_1$ for increased $\xi$. However this happens only for the high quality of additional food. It is observed that the system moves from limit cycle to $E^*$ in the presence of additional food and remains there even for higher values of $\xi$. The parameter values used for this analysis are $\gamma = 17.0,\  \alpha = 10.0,\  \epsilon = 0.35,\  \delta = 3.8,\  m = 3.0$.
	
	From figures \ref{r1xi3} - \ref{r4xi3}, the change in nature of the equilibrium point is observed with increased values of $\xi$. There is also a possibility of bifurcation between $E_1 = (\gamma,0)$ and $E^* = (x^*,y^*)$. 
	
	\newpage
	
	\section{Effect of Mutual Interference on the Dynamics} \label{secmutual}
	
	\begin{figure}[ht]
		\centering
		\includegraphics[width=\textwidth]{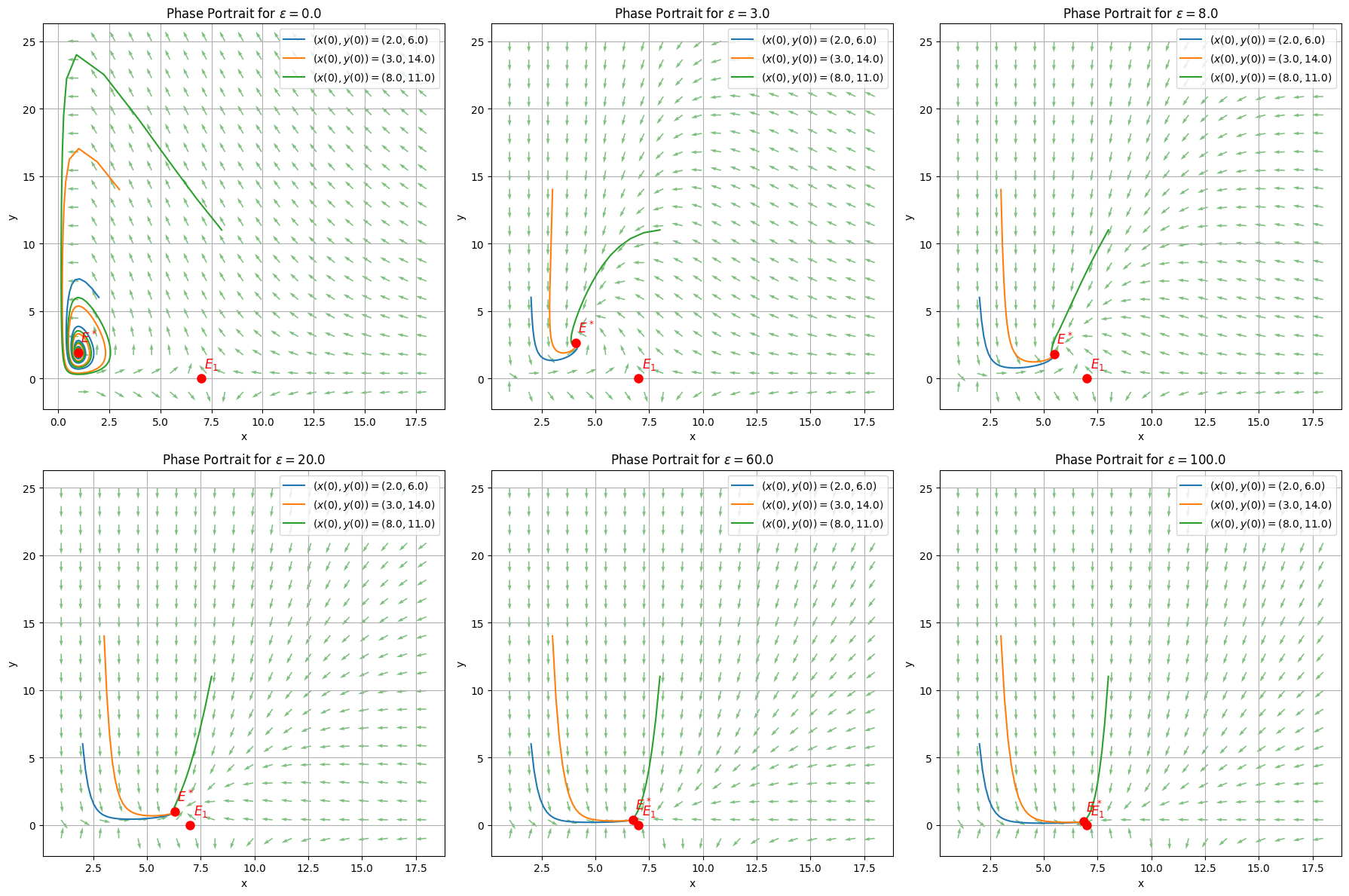}
		\caption{Phase portrait of additional food provided system (\ref{midm3}) for different rates of mutual interference in the presence of low quality and low quantity of additional food.}
		\label{epsll3}
	\end{figure}
	
	\textbf{Low Quality and Low Quantity of additional food:} From figure \ref{epsll3}, it is observed that the system tends to interior equilibrium point even for higher values of $\epsilon$. In the absence of mutual interference, the dynamical system (\ref{midm3}) is driven to $E^*$. The same behavior is observed for different parameter combinations of the system. The parameter values used for this analysis are $\gamma = 7.0, \ \alpha = 0.5,\  \xi = 0.5,\  \delta = 6.0,\  m = 4.0$.
	
	\newpage
	
	\begin{figure}[ht]
		\centering
		\includegraphics[width=\textwidth]{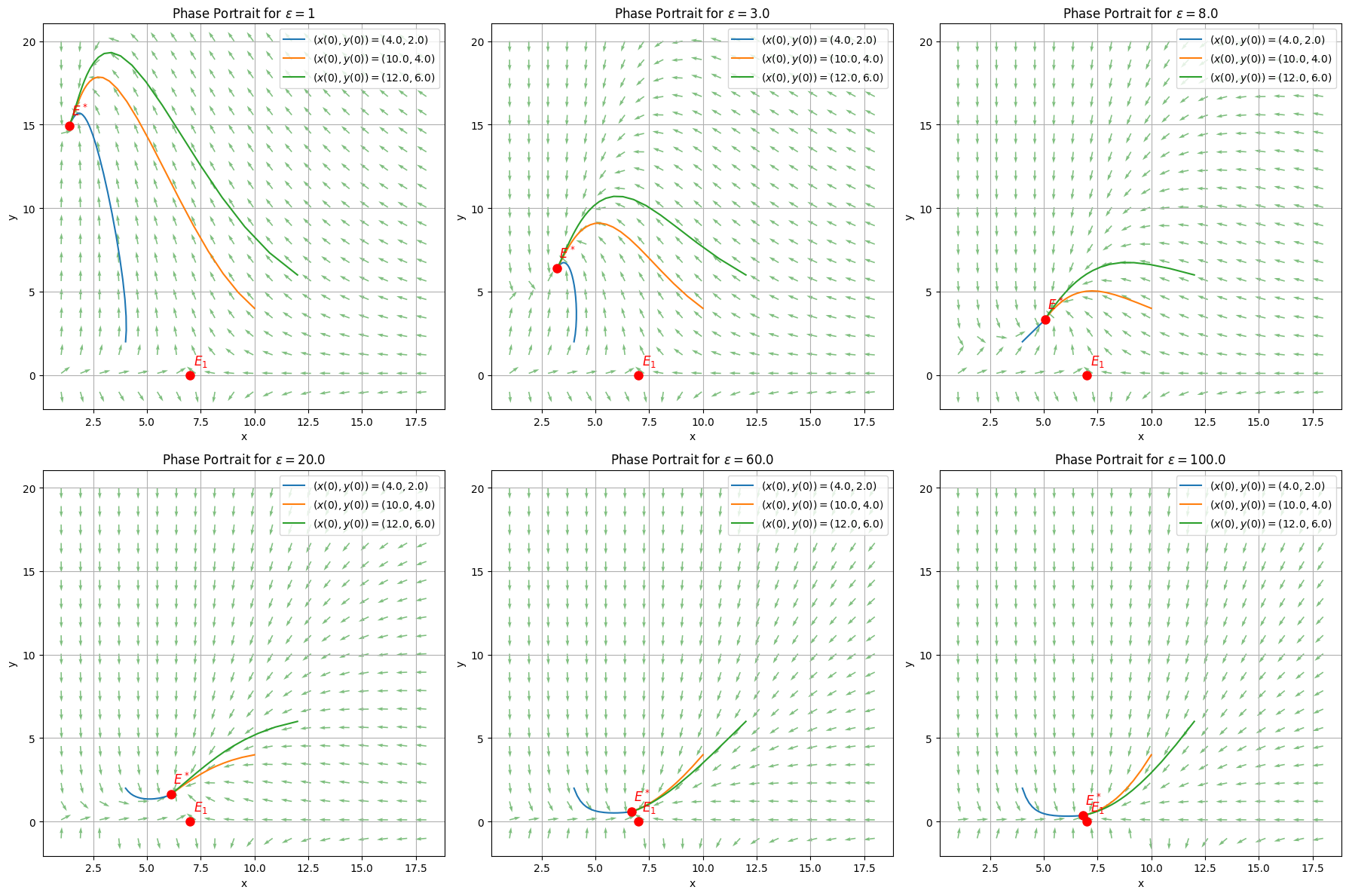}
		\caption{Phase portrait of additional food provided system (\ref{midm3}) for different rates of mutual interference in the presence of low quality and high quantity of additional food.}
		\label{epslh3}
	\end{figure}
	
	\textbf{Low Quality and High Quantity of additional food:} From figure \ref{epslh3}, it is observed that the system tends to interior equilibrium point even for higher values of $\epsilon$. In the absence of mutual interference, the dynamical system (\ref{midm3}) is driven to $E^*$. The same behavior is observed for different parameter combinations of the system. The parameter values used for this analysis are $\gamma = 7.0, \ \alpha = 0.5, \ \xi = 15,\  \delta = 6.0, \ m = 4.0$.
	
	\newpage
	
	\begin{figure}[ht]
		\centering
		\includegraphics[width=\textwidth]{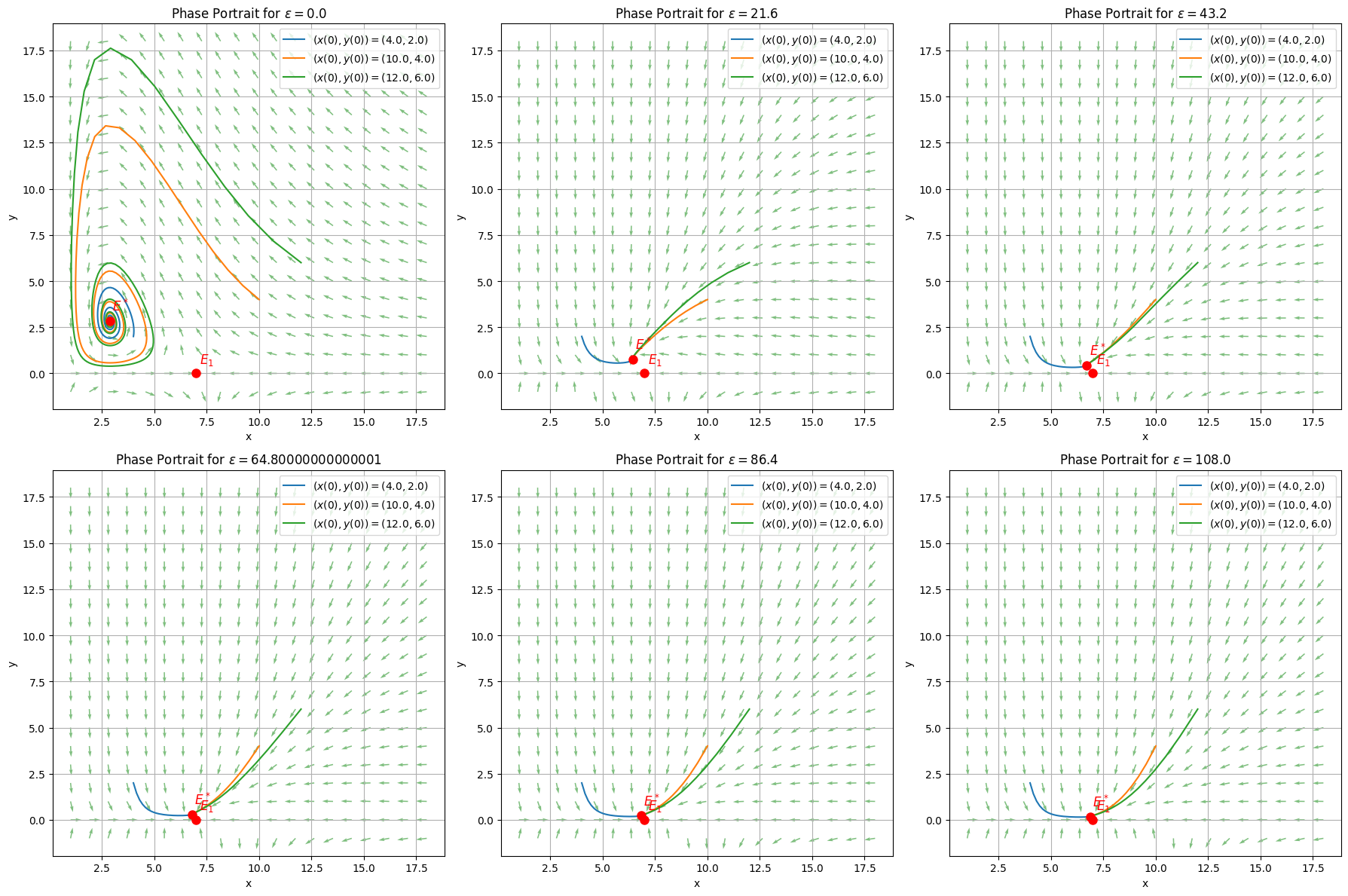}
		\caption{Phase portrait of additional food provided system (\ref{midm3}) for different rates of mutual interference in the presence of high quality and low quantity of additional food.}
		\label{epshl3}
	\end{figure}
	
	\textbf{High Quality and Low Quantity of additional food:} From figure \ref{epshl3}, it is observed that the system tends to interior equilibrium point even for higher values of $\epsilon$. In the absence of mutual interference, the dynamical system (\ref{midm3}) is driven to $E^*$. The same behavior is observed for different parameter combinations of the system. The parameter values used for this analysis are $\gamma = 7.0, \ \alpha = 5.0, \ \xi = 0.9,\  \delta = 6.0,\  m = 4.0$.
	
	\newpage
	
	\begin{figure}[ht]
		\centering
		\includegraphics[width=\textwidth]{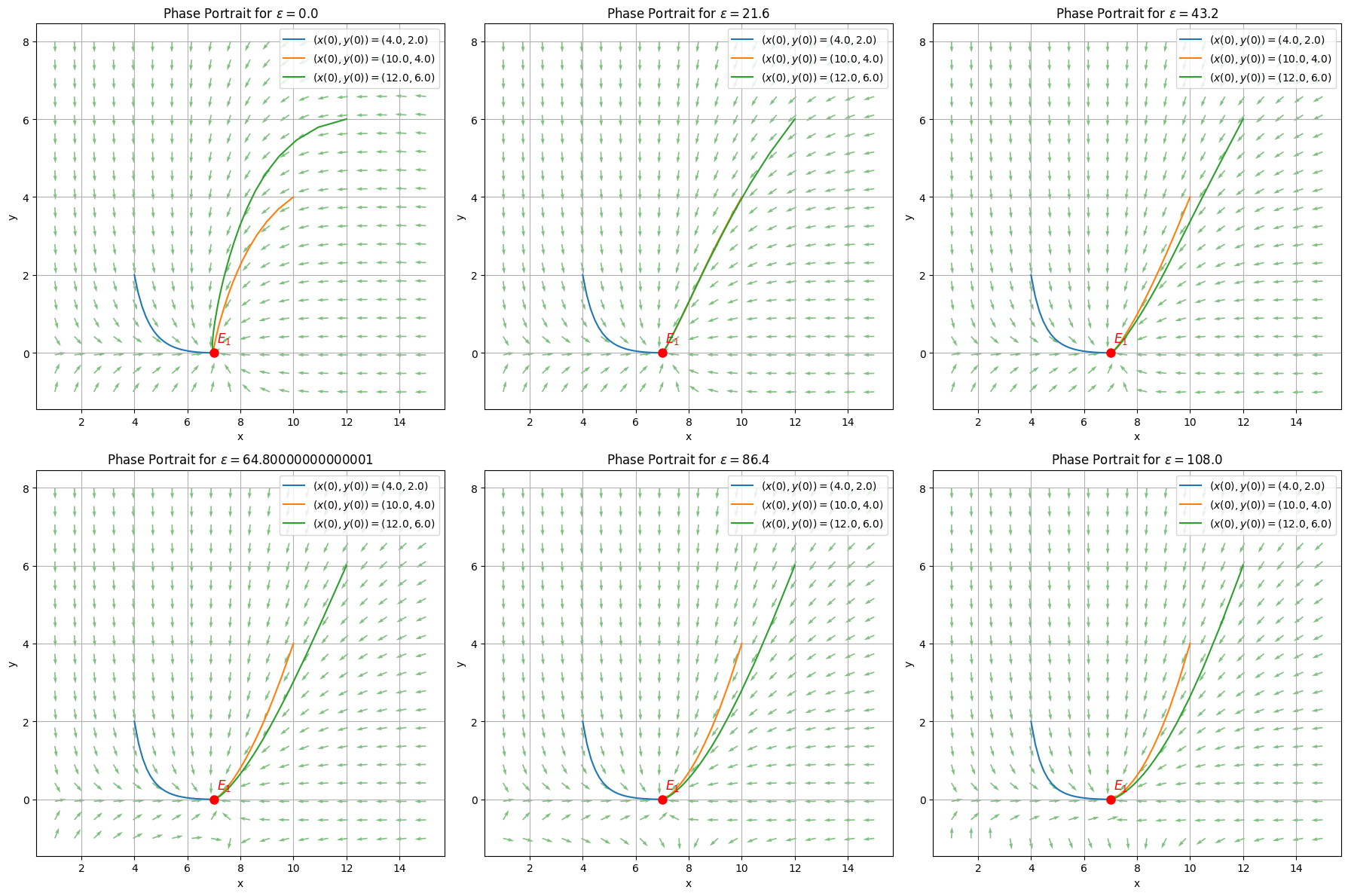}
		\caption{Phase portrait of additional food provided system (\ref{midm3}) for different rates of mutual interference in the presence of high quality and high quantity of additional food.}
		\label{epshh3}
	\end{figure}
	
	\textbf{High Quality and High Quantity of additional food:} From figure \ref{epshh3}, it is observed that the system tends to axial equilibrium point even for higher values of $\epsilon$. In the absence of mutual interference, the dynamical system (\ref{midm3}) is driven to $E_1$. The same behavior is observed for different parameter combinations of the system. The parameter values used for this analysis are $\gamma = 7.0,\  \alpha = 10.0, \ \xi = 10.0,\  \delta = 6.0, \ m = 4.0$.
	
	\newpage
	
	\section{Time-Optimal Control Studies for Holling type-III Systems} \label{seccontrol}
	
	In this section, we formulate and characterise two time-optimal control problems with quality of additional food and quantity of additional food as control parameters respectively. We shall drive the system (\ref{midm3}) from the initial state $(x_0,y_0)$ to the final state $(\bar{x},\bar{y})$ in minimum time.
	
	\subsection{Quality of Additional Food as Control Parameter}
	
	We assume that the quantity of additional food $(\xi)$ is constant and the quality of additional food varies in $[\alpha_{\text{min}},\alpha_{\text{max}}]$. The time-optimal control problem with additional food provided prey-predator system involving Holling type-III functional response among mutually interfering predators (\ref{midm3}) with quality of additional food ($\alpha$) as control parameter is given by
	
	\begin{equation}
		\begin{rcases}
			& \displaystyle {\bf{\min_{\alpha_{\min} \leq \alpha(t) \leq \alpha_{\max}} T}} \\
			& \text{subject to:} \\
			& \dot{x}(t) = x(t) \Bigg( 1 - \frac{x(t)}{\gamma} \Bigg) - \frac{x^2(t) y(t)}{1 + x^2(t) + \alpha (t) \xi + \epsilon y(t)} \\
			& \dot{y}(t) = \delta y(t) \Bigg( \frac{x^2(t) + \xi}{1+x^2(t)+\alpha (t) \xi + \epsilon y(t)} \Bigg) - m y(t) \\
			& (x(0),y(0)) = (x_0,y_0) \ \text{and} \ (x(T),y(T)) = (\bar{x},\bar{y}).
		\end{rcases}
		\label{alpha03}
	\end{equation}
	
	This problem can be solved using a transformation on the independent variable $t$ by introducing an independent variable $s$ such that $\mathrm{d}t = (1 + \alpha \xi + x^2 + \epsilon y) \mathrm{d}s$. This transformation converts the time-optimal control problem ($\ref{alpha03}$) into the following linear problem.
	
	\begin{equation}
		\begin{rcases}
			& \displaystyle {\bf{\min_{\alpha_{\min} \leq \alpha(t) \leq \alpha_{\max}} S}} \\
			& \text{subject to:} \\
			& \mathring{x}(s) = x (1 - \frac{x}{\gamma}) (1 + x^2 + \alpha  \xi + \epsilon y) - x^2 y \\
			& \mathring{y}(s) = \delta (x^2 + \xi) y - (1 + x^2 + \alpha \xi + \epsilon y) m y \\
			& (x(0),y(0)) = (x_0,y_0) \ \text{and} \ (x(S),y(S)) = (\bar{x},\bar{y}).
		\end{rcases}
		\label{alpha3}
	\end{equation}
	
	Hamiltonian function for this problem (\ref{alpha3}) is given by
	\begin{equation*}
		\begin{split}
			\mathbb{H}(s,x,y,p,q) &= p \Big[x \Big(1 - \frac{x}{\gamma}\Big) (1 + x^2 + \alpha  \xi + \epsilon y) - x^2 y \Big] \\ & + q \Big[\delta (x^2 + \xi) y - (1 + x^2 + \alpha \xi + \epsilon y) m y \Big] \\ 
			&= \Big[ p x \Big(1 - \frac{x}{\gamma}\Big) \xi - q \xi m y \Big] \alpha \\ &     + \Big[ p \Big( x \Big(1 - \frac{x}{\gamma}\Big) (1 + x^2+\epsilon y) - x^2 y\Big) + q \Big(\delta (x^2 + \xi) y - (1 + x^2 + \epsilon y) m y \Big)\Big] \\
		\end{split}
	\end{equation*}
	
	Here, $p$ and $q$ are costate variables satisfying the adjoint equations 
	
	\begin{equation*}
		\begin{split}
			\mathring{p} = \frac{\mathrm{d} p}{\mathrm{d} s} = -  \frac{\partial \mathbb{H}}{\partial x} =  & p \Big[2 x y - 2 x^2 \Big(1- \frac{x}{\gamma}\Big) - (1 + x^2 + \alpha \xi + \epsilon y) \Big(1 - \frac{2 x}{\gamma}\Big) \Big] + 2 q x y (m - \delta) \\
			\mathring{q} = \frac{\mathrm{d} q}{\mathrm{d} s} = -  \frac{\partial \mathbb{H}}{\partial y} =  & p \Big[x^2 - \epsilon x \Big(1 - \frac{x}{\gamma} \Big) \Big] + q \Big[ 2 m \epsilon y + m (1+x^2+\alpha \xi) - \delta (x^2 + \xi)  \Big]
		\end{split}
	\end{equation*}
	
	Since Hamiltonian is a linear function in $\alpha$, the existence of optimal control is guarenteed by \textit{Filippov existence theorem} \cite{cesari2012optimization}.  Since we are minimizing the Hamiltonian, the optimal strategy will be a combination of bang-bang and singular controls which is given by 
	
	\begin{equation}
		\alpha^*(t) =
		\begin{cases}
			\alpha_{\max}, &\text{ if } \frac{\partial \mathbb{H}}{\partial \alpha} < 0\\
			\alpha_{\min}, &\text{ if } \frac{\partial \mathbb{H}}{\partial \alpha} > 0
		\end{cases}
	\end{equation}
	where
	\begin{equation}
		\frac{\partial \mathbb{H}}{\partial \alpha} = p x \Big(1 - \frac{x}{\gamma}\Big) \xi - q \xi m y
	\end{equation}
	
	This problem (\ref{alpha3}) admits a singular solution if there exists an interval $[s_1,s_2]$ on which $\frac{\partial \mathbb{H}}{\partial \alpha} = 0$. Therefore, 
	
	\begin{equation}
		\frac{\partial \mathbb{H}}{\partial \alpha} = p x \Big(1 - \frac{x}{\gamma}\Big) \xi - q \xi m y = 0 \textit{ i.e. } \frac{p}{q} = \frac{\gamma m y}{x (\gamma - x)} \label{apbyq13}
	\end{equation}
	
	Differentiating $\frac{\partial \mathbb{H}}{\partial \alpha}$ with respect to $s$ we obtain 
	\begin{equation*}
		\begin{split}    
			\frac{\mathrm{d}}{\mathrm{d}s} \frac{\partial \mathbb{H}}{\partial \alpha} =& \frac{\mathrm{d}}{\mathrm{d}s} \Big[  p x \Big(1 - \frac{x}{\gamma}\Big) \xi - q \xi m y  \Big] \\
			=&  \xi x \Big(1-\frac{x}{\gamma}\Big) \mathring{p} + p \xi \Big(1 - \frac{2 x}{\gamma}\Big) \mathring{x} - \xi m y \mathring{q} - q \xi m \mathring{y}
		\end{split}
	\end{equation*}
	
	Substituting the values of $\mathring{x}, \mathring{y}, \mathring{p}, \mathring{q}$ in the above equation, we obtain
	
	\begin{equation*}
		\begin{split}
			\frac{\mathrm{d}}{\mathrm{d}s} \frac{\partial \mathbb{H}}{\partial \alpha} = & \xi x \Big(1-\frac{x}{\gamma}\Big) \Bigg[ p \Big[2 x y - 2 x^2 \Big(1- \frac{x}{\gamma}\Big) - (1 + x^2 + \alpha \xi + \epsilon y) \Big(1 - \frac{2 x}{\gamma}\Big) \Big] + 2 q x y (m - \delta) \Bigg] \\ & + p \xi \Big(1 - \frac{2 x}{\gamma}\Big) \Bigg[x (1 - \frac{x}{\gamma}) (1 + x^2 + \alpha  \xi + \epsilon y) - x^2 y \Bigg] \\ &  - \xi m y \Bigg[p \Big[x^2 - \epsilon x \Big(1 - \frac{x}{\gamma} \Big) \Big] + q \Big[ 2 m \epsilon y + m (1+x^2+\alpha \xi) - \delta (x^2 + \xi)  \Big] \Bigg] \\ &  - q \xi m \Bigg[ \delta (x^2 + \xi) y - (1 + x^2 + \alpha \xi + \epsilon y) m y \Bigg]
		\end{split}
	\end{equation*}
	
	Along the singular arc, $\frac{\mathrm{d}}{\mathrm{d}s} \frac{\partial \mathbb{H}}{\partial \alpha} = 0$. This implies that 
	
	\begin{equation*}
		\begin{split}
			\xi x \Big(1-\frac{x}{\gamma}\Big) \Bigg[ p \Big[2 x y - 2 x^2 \Big(1- \frac{x}{\gamma}\Big) - (1 + x^2 + \alpha \xi + \epsilon y) \Big(1 - \frac{2 x}{\gamma}\Big) \Big] + 2 q x y (m - \delta) \Bigg] & \\ + p \xi \Big(1 - \frac{2 x}{\gamma}\Big) \Bigg[x \Big(1 - \frac{x}{\gamma}\Big) (1 + x^2 + \alpha  \xi + \epsilon y) - x^2 y \Bigg] & \\ - \xi m y \Bigg[p \Big[x^2 - \epsilon x \Big(1 - \frac{x}{\gamma} \Big) \Big] + q \Big[ 2 m \epsilon y + m (1+x^2+\alpha \xi) - \delta (x^2 + \xi)  \Big] \Bigg] & \\  - q \xi m \Bigg[ \delta (x^2 + \xi) y - (1 + x^2 + \alpha \xi + \epsilon y) m y \Bigg] & = 0
		\end{split}
	\end{equation*}
	
	Upon simplification, we get 
	
	\begin{equation*}
		\begin{split}
			\xi x \Big(1-\frac{x}{\gamma}\Big) \Bigg[ 2 p x \Big[y - x \Big(1- \frac{x}{\gamma}\Big) \Big] + 2 q x y (m - \delta) \Bigg] - p x^2 y \xi \Big(1 - \frac{2 x}{\gamma}\Big)  & \\ - \xi m y \Bigg[p \Big[x^2 - \epsilon x \Big(1 - \frac{x}{\gamma} \Big) \Big] + 2 q m \epsilon y  \Bigg]  + q \xi m^2 \epsilon y^2 & = 0
		\end{split}
	\end{equation*}
	
	and that 
	\begin{equation} \label{apbyq23}
		\frac{p}{q} = \frac{-\xi y \Bigg[ 2 (m - \delta) x^2 \Big(1-\frac{x}{\gamma}\Big) - \epsilon m^2 y \Bigg]}{\xi x \Bigg[ (1-m) x y + \Big(1-\frac{x}{\gamma}\Big) \Big(m \epsilon y - 2 x^2 \Big(1-\frac{x}{\gamma})) \Bigg]}
	\end{equation}
	
	From (\ref{apbyq13}) and (\ref{apbyq23}), we have $$ y = \frac{2 \delta }{m (1 - m)} x \Big( 1 -\frac{x}{\gamma} \Big)^2$$The solutions of this cubic equation gives the switching points of the bang-bang control.
	
	\subsection{Quantity of Additional Food as Control Parameter}
	
	In this section, We assume that the quality of additional food $(\alpha)$ is constant and the quantity of additional food varies in $[\xi_{\text{min}},\xi_{\text{max}}]$. The time-optimal control problem with additional food provided prey-predator system involving Holling type-III functional response among mutually interfering predators (\ref{midm3}) with quantity of additional food ($\xi$) as control parameter is given by
	
	\begin{equation}
		\begin{rcases}
			& \displaystyle {\bf{\min_{\xi_{\min} \leq \xi(t) \leq \xi_{\max}} T}} \\
			& \text{subject to:} \\
			& \dot{x}(t) = x(t) \Bigg( 1 - \frac{x(t)}{\gamma} \Bigg) - \frac{x^2(t) y(t)}{1 + x^2(t) + \alpha (t) \xi + \epsilon y(t)} \\
			& \dot{y}(t) = \delta y(t) \Bigg( \frac{x^2(t) + \xi}{1+x^2(t)+\alpha (t) \xi + \epsilon y(t)} \Bigg) - m y(t) \\
			& (x(0),y(0)) = (x_0,y_0) \ \text{and} \ (x(T),y(T)) = (\bar{x},\bar{y}).
		\end{rcases}
		\label{xi03}
	\end{equation}
	
	This problem can be solved using a transformation on the independent variable $t$ by introducing an independent variable $s$ such that $\mathrm{d}t = (1 + \alpha \xi + x^2) \mathrm{d}s$. This transformation converts the time-optimal control problem ($\ref{xi03}$) into the following linear problem.
	
	\begin{equation}
		\begin{rcases}
			& \displaystyle {\bf{\min_{\xi_{\min} \leq \xi(t) \leq \xi_{\max}} S}} \\
			& \text{subject to:} \\
			& \mathring{x}(s) = x (1 - \frac{x}{\gamma}) (1 + x^2 + \alpha  \xi + \epsilon y) - x^2 y \\
			& \mathring{y}(s) = \delta (x^2 + \xi) y - (1 + x^2 + \alpha \xi + \epsilon y) m y \\
			& (x(0),y(0)) = (x_0,y_0) \ \text{and} \ (x(S),y(S)) = (\bar{x},\bar{y}).
		\end{rcases}
		\label{xi3}
	\end{equation}
	
	Hamiltonian function for this problem (\ref{xi3}) is given by
	\begin{equation*}
		\begin{split}
			\mathbb{H}(s,x,y,p,q) &= p \Big[x \Big(1 - \frac{x}{\gamma}\Big) (1 + x^2 + \alpha  \xi + \epsilon y) - x^2 y \Big] \\ & + q \Big[\delta (x^2 + \xi) y - (1 + x^2 + \alpha \xi + \epsilon y) m y \Big] \\ 
			&= \Big[ p x \Big(1 - \frac{x}{\gamma}\Big) \alpha + q \delta y  - q \alpha m y \Big] \xi \\ &     + \Big[ p \Big( x \Big(1 - \frac{x}{\gamma}\Big) (1 + x^2+\epsilon y) - x^2 y\Big) + q \Big(\delta x^2 y - (1 + x^2 + \epsilon y) m y \Big)\Big] \\
		\end{split}
	\end{equation*}
	
	Here, $p$ and $q$ are costate variables satisfying the adjoint equations 
	
	\begin{equation*}
		\begin{split}
			\mathring{p} = \frac{\mathrm{d} p}{\mathrm{d} s} = -  \frac{\partial \mathbb{H}}{\partial x} =  & p \Big[2 x y - 2 x^2 \Big(1- \frac{x}{\gamma}\Big) - (1 + x^2 + \alpha \xi + \epsilon y) \Big(1 - \frac{2 x}{\gamma}\Big) \Big] + 2 q x y (m - \delta) \\
			\mathring{q} = \frac{\mathrm{d} q}{\mathrm{d} s} = -  \frac{\partial \mathbb{H}}{\partial y} =  & p \Big[x^2 - \epsilon x \Big(1 - \frac{x}{\gamma} \Big) \Big] + q \Big[ 2 m \epsilon y + m (1+x^2+\alpha \xi) - \delta (x^2 + \xi)  \Big]
		\end{split}
	\end{equation*}
	
	Since Hamiltonian is a linear function in $\xi$, the existence of optimal control is guarenteed by \textit{Filippov existence theorem} \cite{cesari2012optimization}.  Since we are minimizing the Hamiltonian, the optimal strategy will be a combination of bang-bang and singular controls which is given by 
	
	\begin{equation}
		\xi^*(t) =
		\begin{cases}
			\xi_{\max}, &\text{ if } \frac{\partial \mathbb{H}}{\partial \xi} < 0\\
			\xi_{\min}, &\text{ if } \frac{\partial \mathbb{H}}{\partial \xi} > 0
		\end{cases}
	\end{equation}
	where
	\begin{equation}
		\frac{\partial \mathbb{H}}{\partial \xi} = \alpha p r x \Big(1 - \frac{x}{\gamma}\Big) + q (g y - \alpha (m y + \delta y^2)) 
	\end{equation}
	
	This problem \ref{xi3} admits a singular solution if there exists an interval $[s_1,s_2]$ on which $\frac{\partial \mathbb{H}}{\partial \xi} = 0$. Therefore, 
	
	\begin{equation}
		\frac{\partial \mathbb{H}}{\partial \xi} = \alpha p r x \Big(1 - \frac{x}{\gamma}\Big) + q (g y - \alpha (m y + \delta y^2)) = 0 \textit{ i.e. } \frac{p}{q} = \frac{\gamma (g y - \alpha m y - \alpha \delta y^2)}{\alpha r x (-\gamma + x)} \label{xpbyq13}
	\end{equation}
	
	Differentiating $\frac{\partial \mathbb{H}}{\partial \xi}$ with respect to $s$ we obtain 
	
	\begin{equation*}
		\begin{split}
			\frac{\mathrm{d}}{\mathrm{d}s} \frac{\partial \mathbb{H}}{\partial \xi} &= \frac{\mathrm{d}}{\mathrm{d}s} \Big[ \alpha p r x \Big(1 - \frac{x}{\gamma}\Big) + q \Big(g y - \alpha (m y + \delta y^2)\Big) \Big] \\
			&= \alpha p r \Big(1-\frac{2 x}{\gamma}\Big) \dot{x} + q (g -\alpha (m +\delta 2 y )) \dot{y} + \alpha r x \Big(1-\frac{x}{\gamma}\Big) \dot{p} + (g y-\alpha (m y+\delta y^2)) \dot{q}
		\end{split}
	\end{equation*}
	
	Substituting the values of $\dot{x}, \dot{y}, \dot{p}, \dot{q}$ in the above equation, we obtain
	\begin{equation*}
		\begin{split}
			\frac{\mathrm{d}}{\mathrm{d}s} \frac{\partial \mathbb{H}}{\partial \xi} =& \alpha p r \Big(1 - \frac{2 x}{\gamma}\Big) \Big(r x \Big(1 - \frac{x}{\gamma}\Big) (1 + x^2 + \alpha \xi) - x^2 y \Big) \\ & + \Big(g y - \alpha (m y + \delta y^2)\Big) \Big(p x^2 - q \Big(g (x^2 + \xi) - (1 + x^2 + \alpha \xi) (m + 2 \delta y)\Big)\Big) \\ & +  \alpha r x \Big(1 - \frac{x}{\gamma}\Big) \Big[-p (2 r x^2 \Big(1 - \frac{x}{\gamma}\Big) - \frac{r x (1 + x^2 + \alpha \xi)}{\gamma} + r (1 - \frac{x}{\gamma}) (1 + x^2 + \alpha \xi) - 2 x y) \\ & - q (2 g x y - 2 x (m y + \delta y^2))\Big] \\ &  + q (g (g (x^2 + \xi) y - (1 + x^2 + \alpha \xi) (m y + \delta y^2)) - \alpha (m (g (x^2 + \xi) y - (1 + x^2 + \alpha \xi) (m y + \delta y^2)) \\ & + 2 \delta y (g (x^2 + \xi) y - (1 + x^2 + \alpha \xi) (m y + \delta y^2))))
		\end{split}
	\end{equation*}
	
	Along the singular arc, $\frac{\mathrm{d}}{\mathrm{d}s} \frac{\partial \mathbb{H}}{\partial \xi} = 0$. This implies that 
	
	\begin{equation*}
		\begin{split}
			\alpha p r \Big(1 - \frac{2 x}{\gamma}\Big) \Big(r x \Big(1 - \frac{x}{\gamma}\Big) (1 + x^2 + \alpha \xi) - x^2 y \Big) & \\ + \Big(g y - \alpha (m y + \delta y^2)\Big) \Big(p x^2 - q \Big(g (x^2 + \xi) - (1 + x^2 + \alpha \xi) (m + 2 \delta y)\Big)\Big) & \\ +  \alpha r x \Big(1 - \frac{x}{\gamma}\Big) \Big[-p (2 r x^2 \Big(1 - \frac{x}{\gamma}\Big) - \frac{r x (1 + x^2 + \alpha \xi)}{\gamma} + r (1 - \frac{x}{\gamma}) (1 + x^2 + \alpha \xi) - 2 x y) &  \\ - q (2 g x y - 2 x (m y + \delta y^2))\Big] & \\  + q (g (g (x^2 + \xi) y - (1 + x^2 + \alpha \xi) (m y + \delta y^2)) - \alpha (m (g (x^2 + \xi) y - (1 + x^2 + \alpha \xi) (m y + \delta y^2)) & \\  + 2 \delta y (g (x^2 + \xi) y - (1 + x^2 + \alpha \xi) (m y + \delta y^2)))) & = 0
		\end{split}
	\end{equation*}
	
	and that 
	\begin{equation} \label{xpbyq23}
		\frac{p}{q} = \frac{\gamma y (\delta g \gamma (1 + x^2) y + \alpha x^2 (2 r (\gamma - x) (m + \delta y) - g (2 \gamma r - 2 r x + \delta \gamma y)))}{x^2 (2 \alpha r^2 (\gamma - x)^2 x - \gamma^2 (g + \alpha (-m + r)) y + \alpha \delta \gamma^2 y^2)}
	\end{equation}
	
	From (\ref{xpbyq13}) and (\ref{xpbyq23}), we have $$\frac{x y (-g + \alpha (m + \delta y))}{\alpha r (\gamma - x)} + \frac{y (-\delta g \gamma (1 + x^2) y + \alpha x^2 (-2 r (\gamma - x) (m + \delta y) + g (2 \gamma r - 2 r x + \delta \gamma y)))}{2 \alpha r^2 (\gamma - x)^2 x - \gamma^2 (g + \alpha (-m + r)) y + \alpha \delta \gamma^2 y^2} = 0.$$ The solutions of this cubic equation gives the switching points of the bang-bang control.
	
	\section{Numerical Simulations} \label{secnumer}
	
	In this section, we will numerically illustrate the results of theoretical findings obtained for the time optimal control problem. We use the object oriented Python library, GEKKO optimization suite \cite{beal2018gekko}, to solve the time optimal control problems (\ref{alpha03}) and (\ref{xi03}) numerically. We use the non-simulation simultaneous method provided by GEKKO library to optimize the objective and implicitly calculate the state variables.
	
	\begin{figure}[ht]
		\centering
		\includegraphics[width=\textwidth]{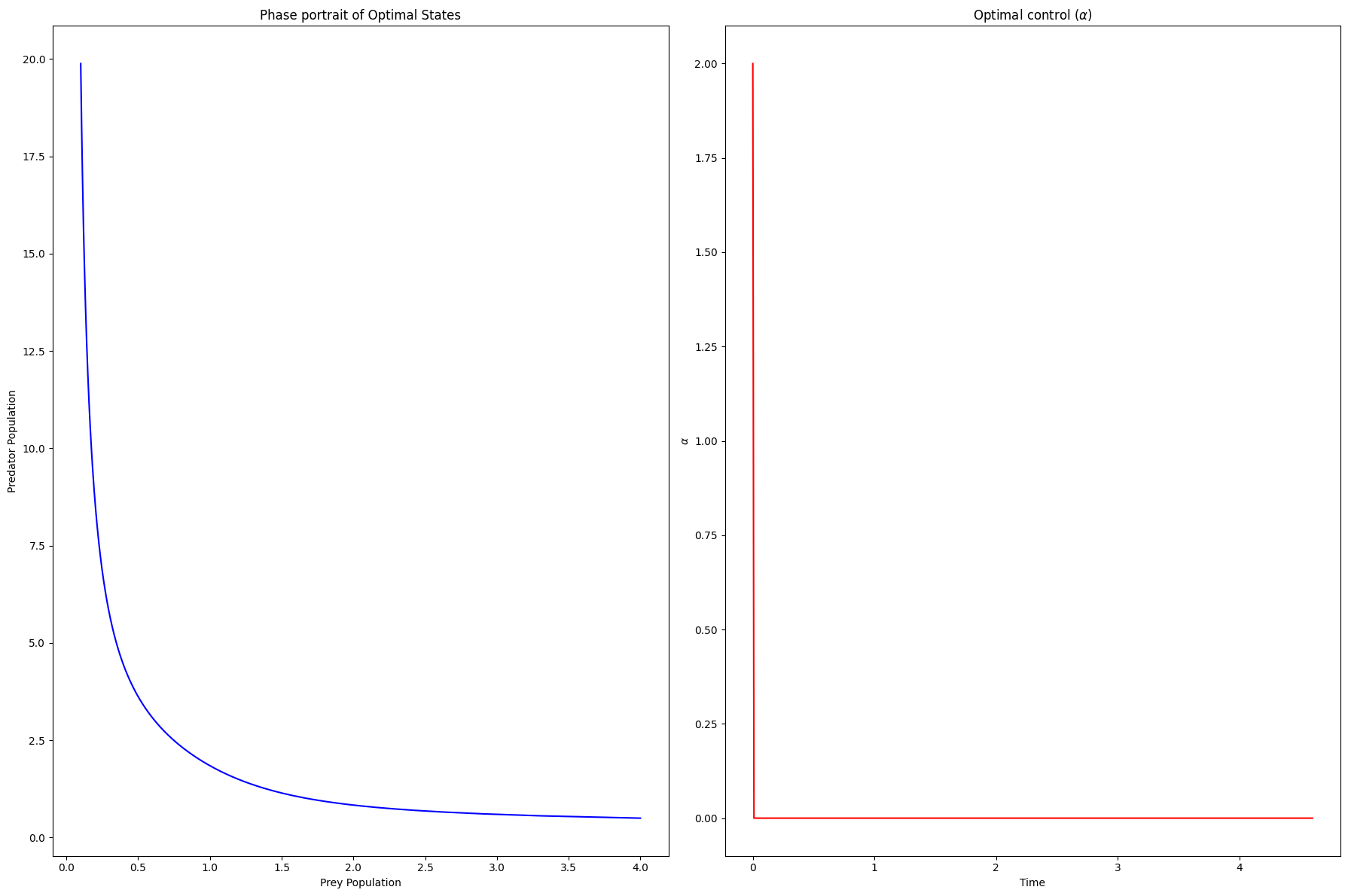}
		\caption{Phase portrait of optimal state vectors and optimal control ($\alpha$) for the time optimal control problem (\ref{alpha03}).}
		\label{alphacontrol3}
	\end{figure}
	
\begin{figure}[ht]
		\centering
		\includegraphics[width=\textwidth]{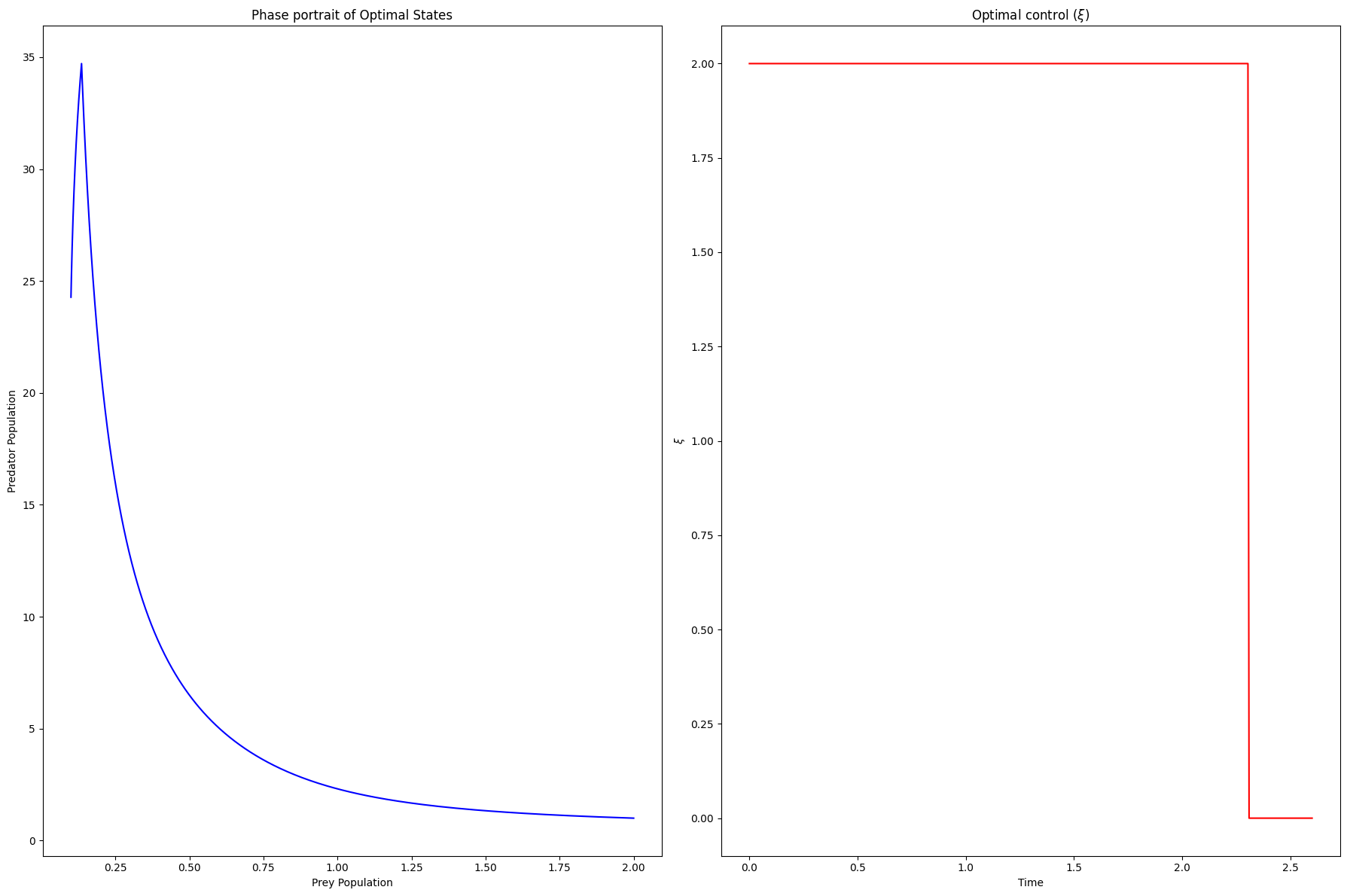}
		\caption{Phase portrait of optimal state vectors and optimal control ($\xi$) for the time optimal control problem (\ref{xi03}).}
		\label{xicontrol3}
	\end{figure}

	Figure \ref{alphacontrol3} represents the optimal state and optimal control trajectories of the time optimal control problem (\ref{alpha03}). In this simulation, control variable ($\alpha$) ranges from $0$ to $2$. For the parameter values $\gamma = 1.2,\  \epsilon = 0.05,\  \xi = 0.8,\  \delta = 1.5,\  m = 0.5$, the optimal control $\alpha$ steers the system from the initial point $(4.0,0.5)$ to final point $(0.1,19.7)$. From the figure \ref{xicontrol3}, it is observed that the optimal control is a bang-bang control. Here, the bang-bang control drives system to nearly prey-elimination stage. 
	
	Figure \ref{xicontrol3} represents the optimal state and optimal control trajectories of the time optimal control problem (\ref{xi03}). In this simulation, control variable ($\xi$) ranges from $0$ to $2$. For the parameter values $\gamma = 1.2, \ \epsilon = 0.05,\  \alpha = 1.0,\  \delta = 1.5,\  m = 0.5$, the optimal control $\xi$ steers the system from the initial point $(2.0,1.0)$ to final point $(0.1,24.5)$. From the figure \ref{xicontrol3}, it is observed that the optimal control is a bang-bang control. Here, the bang-bang control drives system to nearly prey-elimination stage. These results validate the theoretical findings.

	\newpage
	
	\section{Discussions and Conclusions} \label{secdisc}
	
	\qquad In this paper, we derived the additional food provided prey-predator system exhibiting Holling type-III functional response among mutually interfering predators. From the positivity and boundedness result, we proved that this system remains in the positive-$xy$ quadrant when the system starts from the positive-$xy$ quadrant. Since we are working on the population of pest and natural enemies, populations being positive and bounded proves the fundamental assumption on populations. Later, we performed the existence and stability of equilibria of both the initial system and additional food provided system. In this analysis, we not only listed the possible equilibrium points in the positive-$xy$ quadrant but also discussed the stability of these multiple equilibrium points. To sum up all the results of stability analysis, we showed that the additional food provided system has interior equilibrium point as a stable equilibrium when $\xi > \frac{m - (\delta - m) \gamma^2}{\delta - m \alpha}$. In the other case, it has axial equilibrium $(\gamma,0)$ as it's stable equilibrium. 
	
	Further, we explored the global dynamics and stability behavior exhibited by the prey-predator model. The parameter space of the initial system is divided into four regions, where the initial system ($\ref{midm03}$) tends to axial equilibrium, interior equilibrium, interior equilibrium and stable limit cycle respectively. We explored the influence of additional food parameters $\alpha$ and $\xi$ in each of these regions numerically. We could show the impact of additional food on shifting the equilibrium point from that of initial system in each of the regions. We also performed a numerical study on the influence of mutual interference term in four possible cases of low and high, quality and quantity of additional food respectively. 
	
	In the end, we framed two time optimal control problems where the additional food provided system is the state space and the quality of additional food and the quantity of additional food as control parameters respectively. In order to implement the biological control of pests, we chose prey as pests and predators as natural enemies. It is observed that pests reach nearly zero quantity in finite time and the optimal control is a bang-bang control. The obtained numerical results are also in line with the theoretical findings. We proved that the bang-bang control is the optimal control for achieving pest elimination through biological control.
	
	Taking the analysis even further, we wish to explore the possibility of occurrence of multi dimensional bifurcations in our future works. This in-depth mathematical analysis enhances our understanding of the dynamical system. In addition to this, studies on the impact of simultaneously applied multiple controls can delve into the better strategies for biological control of pest populations. It is also observed that the system parameters exhibit random fluctuations due to environment and the results obtained on these advanced stochastic models will be a valuable addition to the existing knowledge of mathematical modeling for biological control of pests.

	\paragraph{Author contributions}D.B.P. contributed to the conceptualization, conducted all mathematical and numerical analyses, generated the figures and drafted the manuscript. D.K.K.V. contributed to the conceptualization, research analysis and edited  the flow of manuscript. All authors reviewed and approved the final manuscript. 
	
	\paragraph{Funding}This research was funded by National Board of Higher Mathematics (NBHM), Government of India under project grant - {\bf{Time Optimal Control and Bifurcation Analysis of Coupled Nonlinear Dynamical Systems with Applications to Pest Management, \\ Sanction number: (02011/11/2021NBHM(R.P)/R$\&$D II/10074).}} The funder had no role in study design, data collection and analysis, decision to publish or preparation of the manuscript.
	
	\paragraph{Data availability}We do not analyse or generate any datasets, because our work proceeds within a theoretical and mathematical approach.
	
	\paragraph{Acknowledgements}The authors dedicate this paper to the founder chancellor of SSSIHL, Bhagawan Sri Sathya Sai Baba. The corresponding author also dedicates this paper to his loving elder brother D. A. C. Prakash who still lives in his heart.
	
	\section*{Declarations}
	
	\paragraph{Ethics approval and consent to participate}Not applicable.
	
	\paragraph{Consent for publication}Not applicable.
	
	\paragraph{Competing interests}The authors declare no competing interests.
	
	\printbibliography
\end{document}